\documentclass[11pt]{article}
\usepackage{amsmath}
\usepackage{amssymb}
\usepackage{amsthm}
\usepackage{amscd}
\usepackage{amsfonts}
\usepackage{dsfont}
\usepackage{graphicx}%
\usepackage{subcaption}
\usepackage{fancyhdr}
\usepackage{mathtools}
\DeclarePairedDelimiter\ceil{\lceil}{\rceil}
\DeclarePairedDelimiter\floor{\lfloor}{\rfloor}
\usepackage{cite}
\usepackage{xcolor}
\usepackage[ruled]{algorithm2e}
\usepackage{ifthen}
\usepackage{algorithmic}
\newtheorem{theorem}{Theorem}
\newtheorem{lemma}{Lemma}
\newtheorem{remark}{Remark}
\newtheorem{proposition}{Proposition}
\newtheorem{corollary}{Corollary}

\newtheorem{definition}{Definition}
\newtheorem{assumption}{Assumption}

\usepackage{empheq}
\def\scr#1{{\cal #1}} 
\ifodd 1
\def\eq#1{\begin{equation}#1\end{equation}}
\newcommand{\R}{{\rm I\!R}}
\newcommand{\bbb}{\mathbb}

\else

\newcommand{\com}[1]{}
\newcommand{\clar}[1]{}
\newcommand{\response}[1]{}
\fi

\newcommand{\1}{\mathbf{1}}
\newcommand{\0}{\mathbf{0}}

\usepackage{hyperref}

\def\##1\#{\begin{align}#1\end{align}}
\def\$#1\${\begin{align*}#1\end{align*}}

\def\qed{ \rule{.08in}{.08in}}

\newcommand{\dfb}{\stackrel{\Delta}{=}}

\usepackage{parskip}
\title{An Analysis Tool for Push-Sum Based Distributed Optimization 
}

\author{Yixuan Lin \hspace{.3in} Ji Liu\thanks{
Y.~Lin is with the Department of Applied Mathematics and Statistics at Stony Brook University (\texttt{yixuan.liu.1@stonybrook.edu}).
J.~Liu is with the Department of Electrical and Computer Engineering at Stony Brook University
(\texttt{ji.liu@stonybrook.edu}).
}
}

\begin{document}

\maketitle
\thispagestyle{empty}

\begin{abstract}

The push-sum algorithm is probably the most important distributed averaging approach over directed graphs, which has been applied to various problems including distributed optimization. This paper establishes the explicit absolute probability sequence for the push-sum algorithm, and based on which, constructs quadratic Lyapunov functions for push-sum based distributed optimization algorithms. As illustrative examples, the proposed novel analysis tool can improve the convergence rates of the subgradient-push and stochastic gradient-push, two important algorithms for distributed convex optimization over unbalanced directed graphs. Specifically, the paper proves that the subgradient-push algorithm converges at a rate of $O(1/\sqrt{t})$ for general convex functions and stochastic gradient-push algorithm converges at a rate of $O(1/t)$ for strongly convex functions, over time-varying unbalanced directed graphs. Both rates are respectively the same as the state-of-the-art rates of their single-agent counterparts and thus optimal, which closes the theoretical gap between the centralized and push-sum based (sub)gradient methods. The paper further proposes a heterogeneous push-sum based subgradient algorithm in which each agent can arbitrarily switch between subgradient-push and push-subgradient. The heterogeneous algorithm thus subsumes both subgradient-push and push-subgradient as special cases, and still converges to an optimal point at an optimal rate. The proposed tool can also be extended to analyze distributed weighted averaging. 
\end{abstract}

\section{Introduction}

There are three major information fusion schemes in the vast distributed algorithms literature: consensus via stochastic matrices \cite{flocking}, distributed averaging via doubly stochastic matrices \cite{fast}, and push-sum via column stochastic matrices \cite{push}.\footnote{A square nonnegative matrix is called a row stochastic matrix, or simply stochastic matrix, if its row sums all equal one. Similarly, a square nonnegative matrix is called a column stochastic matrix if its column sums all equal one.
A square nonnegative matrix is called a doubly stochastic matrix if its row sums and column sums all equal one.} Among the three, the push-sum scheme is the only one that is able to not only achieve agreement on the average, but also works for directed graphs, allowing uni-directional communication. Because of this, the push-sum scheme has been widely utilized in various distributed algorithms including distributed optimization \cite{tsianos,nedic,nedic2016stochastic}, distributed deep learning \cite{assran19a}, and distributed reinforcement learning \cite{yixuan,jingxuan,finite}.

The push-sum algorithm was first proposed in \cite{push} and later developed in \cite{weighted,wallerton,acc12}, because of which the algorithm is 
sometimes also called weighted gossip \cite{weighted}, ratio consensus \cite{wallerton}, and double linear iterations \cite{acc12}. Push-sum is the most popular and probably the most powerful existing approach to design distributed algorithms over time-varying directed graphs, which can be tailored to be robust against communication delay \cite{hadjicostis2013average}, asynchronous updating \cite{wallerton}, and package drops \cite{robustpush}.
Although the analysis of the push-sum algorithm is elegant, the analyses of push-sum based algorithms are often quite complicated, e.g., subgradient-push \cite{nedic}, DEXTRA \cite{xi2017dextra} (a push-sum based variant of the well-known EXTRA algorithm \cite{shi2015extra}) and Push-DIGing \cite{nedic2017achieving}. Actually, all these push-sum based algorithms rely on the pioneering analysis and results in \cite{nedic}.

Distributed optimization originated from the work of \cite{nedic2009distributed} and has
achieved great success in both theory and practice, as evidenced in survey papers \cite{yang2019survey,nedic2018distributed,molzahn2017survey}. The efforts have been made to design distributed multi-agent versions for various optimization algorithms, including the subgradient method \cite{nedic2009distributed}, alternating direction method of multipliers (ADMM) \cite{wei2012distributed}, Nesterov accelerated gradient method \cite{qu2019accelerated}, and proximal gradient descent \cite{8752033}, to name a few. 
Most existing distributed optimization algorithms require the underlying communication network be described by a possibly time-dependent undirected graph or balanced directed graph\footnote{A weighted directed graph is called balanced if the sum of all in-weights equals the sum of all out-weights at each of its vertices \cite{gharesifard2013distributed}.}, which allows a distributed manner to construct a doubly stochastic matrix \cite{metro2,weightbalance}. Such a distributed algorithm usually achieves the same order of convergence rate as its single-agent counterpart, only with a difference at a constant coefficient which depends on graph connectivity \cite{nedic2018network}.

The push-sum based subgradient algorithm proposed by Nedi\'c and Olshevsky in \cite{nedic} is the first distributed convex optimization algorithm which works for unbalanced directed graphs. There are two ``gaps'' in the analysis in \cite{nedic}. First, the convergence rate analysis is based on a special convex combination of the history of the states of all agents (see Theorem 2 in \cite{nedic}), which is ``unusual'' compared with non-push-sum based distributed optimization algorithms (see e.g. \cite{nedic2009distributed}). Second, more importantly, the convergence rate derived in \cite{nedic} is of order $O(\ln t/\sqrt{t})$, which is slower than that of the single-agent subgradient method, $O(1/\sqrt{t})$ (see e.g. Theorem 7 in \cite{nedic2018network}). In the subsequent work \cite{nedic2016stochastic}, Nedi\'c and Olshevsky further studied a stochastic gradient-push algorithm and showed an $O(\ln t/t)$ convergence rate for strongly convex functions by again studying an untypical convex combination of agents' states.
Apparently there is a convergence rate gap between multi-agent stochastic gradient-push over time-varying directed graphs and state-of-the-art single-agent stochastic gradient descent (see e.g. $O(1/t)$ in \cite{rakhlin2012making}), which is also stated as an open problem in \cite{nedic2016stochastic}.
The problem was tackled in \cite{alex_jmlr} by carefully choosing stepsizes, but with a global quantity available to each agent. 
The preceding discussion leads to the following question:

\begin{center}
    \begin{minipage}{0.95\linewidth}
    \centering
    {\em Can push-sum based distributed subgradient achieve the best possible convergence rate for time-varying, unbalanced, directed graphs, without any network-wide information?}
    \end{minipage}
\end{center}

\vspace{.05in}

This is precisely what we consider in this paper. Indeed, we theoretically prove that the subgradient-push algorithm proposed in \cite{nedic} converges at a rate of $O(1/\sqrt{t})$ for general convex functions, and the stochastic gradient-push algorithm studied in \cite{nedic2016stochastic} converges at a rate of $O(1/t)$ for strongly convex functions, both over time-varying unbalanced directed graphs with appropriate joint connectivity, by analyzing the ``standard'' convex combination of the history of the states of all agents. 
Both rates are respectively the same as the state-of-the-art rates of their single-agent counterparts \cite{nedic2018network,rakhlin2012making} and thus optimal, which closes the theoretical gap between the centralized single-agent and push-sum based multi-agent (sub)gradient methods.

To this end, we first establish the explicit ``absolute probability sequence'' for the push-sum algorithm, and then based on which, construct a quadratic Lyapunov function for the push-sum based subgradient algorithms. The approach yields a novel analysis tool for push-sum based
distributed first-order optimization algorithms over possibly time-varying, unbalanced, directed graphs, and is expected to be applicable to many other push-sum based distributed optimization, computation, and machine learning algorithms.
The ``absolute probability sequence'' concept has been applied to analyze consensus processes \cite{touri2012product,tacrate} and consensus-based distributed optimization \cite{saadatniaki2020decentralized}. To our knowledge, the concept has never been used to analyze the push-sum algorithm and push-sum based distributed algorithms.

\section{Push-Sum}\label{sec:pushsum}

We begin with the push-sum algorithm.
Consider a network consisting of $n$ agents, labeled $1$ through $n$ for the purpose of presentation. The agents are not aware of such a global labeling, but can differentiate between their neighbors. The neighbor relations among the $n$ agents are characterized by a time-dependent directed graph $\bbb{G}(t) = (\mathcal{V},\mathcal{E}(t))$ whose
vertices correspond to agents and whose directed edges (or arcs) depict neighbor relations, where $\mathcal{V}=\{1,\ldots,n\}$ is the vertex set and $\mathcal{E}(t)\subset\mathcal{V} \times \mathcal{V}$ is the directed edge set at time $t$.
Specifically, agent $j$ is an in-neighbor of agent $i$ at time $t$ if $(j,i)\in\scr{E}(t)$, and similarly, agent $k$ is an out-neighbor of agent $i$ at time $t$ if $(i,k)\in\scr{E}(t)$.
Each agent can send information to its out-neighbors and receive information from its in-neighbors. Thus, the directions of edges represent the
directions of information flow. For convenience, we assume that each agent is always an in- and out-neighbor of itself, which implies that $\bbb{G}(t)$ has self-arcs at all vertices for all time $t$. We use $\mathcal{N}_i(t)$ and $\mathcal{N}_i^{-}(t)$ to denote the in- and out-neighbor set of agent $i$ at time $t$, respectively, i.e.,
\begin{align*}
    \mathcal{N}_i(t) = \{ j \in \mathcal{V} \; : \;( j, i ) \in \mathcal{E}(t) \}, \;\;\;\;\;
    \mathcal{N}_i^{-}(t) = \{ k \in \mathcal{V} \; : \; ( i, k ) \in \mathcal{E}(t) \}.
\end{align*}
It is clear that $\mathcal{N}_i(t)$ and $\mathcal{N}_i^{-}(t)$ are nonempty as they both contain index $i$.
In the push-sum algorithm, each agent $i$ has control over two variables, $x_i(t)\in\R^d$ and $y_i(t)\in\R$. At each time $t\in\{0,1,2,\ldots\}$, each agent $j$ transmits two pieces of information, $w_{ij}(t)x_j(t)$ and $w_{ij}(t)y_j(t)$, to its out-neighbour $i$, and then each agent $i$ updates its variables as follows:  
\begin{align}
    x_i(t+1) &= \sum_{j\in\scr{N}_i(t)} w_{ij}(t)x_j(t),\;\;\;\;\; x_i(0)\in\R^d,\label{pushsumx}\\ 
    y_i(t+1) &= \sum_{j\in\scr{N}_i(t)} w_{ij}(t)y_j(t),\;\;\;\;\; y_i(0)=1, \label{pushsumy}
\end{align}
where $w_{ij}(t)$, $j\in\scr{N}_i(t)$, are positive weights satisfying the following assumption.  

\begin{assumption}\label{assum:weighted matrix}
There exists a constant $\beta>0$ such that for all $i,j\in\scr V$ and $t$, $w_{ij}(t) \ge \beta$ whenever $j\in\scr{N}_i(t)$. For all $i\in\scr V$ and $t$, $\sum_{j\in\scr{N}_i^{-}(t)} w_{ji}(t) = 1$.
\end{assumption}

A typical choice of $w_{ij}(t)$ is $1/|\mathcal{N}_j^{-}(t)|$ for all $j\in\scr{N}_i(t)$ which can be computed in a distributed manner and satisfies Assumption~\ref{assum:weighted matrix} with $\beta=1/n$.

Let $W(t)$ be the $n\times n$ matrix whose $ij$th entry equals $w_{ij}(t)$ if $j\in\scr{N}_i(t)$ and zero otherwise; in other words, we set $w_{ij}(t)=0$ for all $j\notin\scr{N}_i(t)$.
From Assumption~\ref{assum:weighted matrix}, each $W(t)$ is a column stochastic matrix that is compliant with the neighbor graph $\bbb{G}(t)$. 
Since each agent $i$ is always assumed to be an in-neighbor of itself, all diagonal entries of $W(t)$ are positive. Thus, if $\bbb{G}(t)$ is strongly connected\footnote{A directed graph is strongly connected if it has a directed path from any vertex to any other vertex.}, $W(t)$ is  irreducible and aperiodic. 
%
Let 
\begin{align*}
    x(t)\dfb \begin{bmatrix}x_1^\top(t)\cr\vdots\cr x_{n}^\top(t)\end{bmatrix}\in\R^{n\times d}, \;\;\;\;\; y(t)\dfb \begin{bmatrix}y_1(t)\cr\vdots\cr y_{n}(t)\end{bmatrix}\in\R^n.
\end{align*}
From \eqref{pushsumx} and \eqref{pushsumy}, $x(t+1)=W(t)x(t)$ and $y(t+1)=W(t)y(t)$. Since $W(t)$ is always a column stochastic matrix for all $t\ge 0$, it is easy to show that $\sum_{i=1}^n x_i(t)=\sum_{i=1}^n x_i(0)$ and $\sum_{i=1}^n y_i(t)=\sum_{i=1}^n y_i(0)=n$ for all $t\ge 0$.

To state the convergence result of the push-sum algorithm \eqref{pushsumx}--\eqref{pushsumy}, we need the following concept.

\begin{definition}\label{def:uniformly}
     A directed graph sequence $\{ \bbb{G}(t) \}$ is uniformly strongly connected if there exists a positive integer $L$ such that for any $t\ge 0$, the union graph $\cup_{k=t}^{t+L-1} \bbb{G}(k)$ is strongly connected.\footnote{The union of two directed graphs, $\bbb G_p$ and $\bbb G_q$, with the same vertex set, written $\bbb G_p \cup \bbb G_q$, is meant the directed graph with the same vertex set and edge set being the union of the edge set of $\bbb G_p$ and $\bbb G_q$. Since this union is a commutative and associative binary operation, the definition extends unambiguously to any finite sequence of directed graphs with the same vertex set.}
     If such an integer exists, we sometimes say that $\{ \bbb{G}(t) \}$ is uniformly strongly connected by sub-sequences of length $L$.
\end{definition}

It is not hard to prove that the above definition is equivalent to the two popular joint connectivity definitions in consensus literature, namely ``$B$-connected'' \cite{nedic2009distributed_quan} and ``repeatedly jointly strongly connected'' \cite{reachingp1}.


\begin{theorem}\label{thm:pushsum}
    If $\{ \bbb{G}(t) \}$ is uniformly strongly connected, then $x_i(t)/y_i(t)$ converges to $\frac{1}{n}\sum_{i=1}^n x_i(0)$ exponentially fast. 
\end{theorem}

The theorem is more or less well-known and its proof can be found e.g. in \cite{push,weighted,wallerton}. 
Here we provide an alternative proof.

We begin with the following weak convergence property of the backward product of an infinite sequence of column stochastic matrices.


\begin{lemma} \label{lemma:pushsum_product}
    If $\{ \bbb{G}(t) \}$ is uniformly strongly connected, then for any fixed $\tau\ge 0$, 
    $W(t)\cdots W(\tau+1)W(\tau)$
    will converge to the set $\{v\1^\top \; :\; v\in\R^n, \1^\top v=1, v>\0\}$ exponentially fast as $t\rightarrow\infty$.\footnote{We use $\0$ and $\1$ to denote the column vectors whose entries all equal to $0$ or $1$, respectively, where the dimensions of the vectors are to be understood from the context. We use $v>\0$ to denote a positive column vector $v$, i.e., each entry of $v$ is positive.}
\end{lemma}

The lemma is essentially the same as Corollary~2~(a) in~\cite{nedic}. Suppose that $\{ \bbb{G}(t) \}$ is uniformly strongly connected by sub-sequences of length $L$. Lemma \ref{lemma:pushsum_product}  implies that there exist constants $c>0$ and $\mu \in [0,1)$ and a sequence of stochastic vectors\footnote{A vector is called a stochastic vector if
its entries are all nonnegative and sum to one.} $ \{ v(t)\}$ such that for all $i,j \in \mathcal{V}$ and $t \ge \tau \ge 0$,
\begin{align}\label{mu}
    \big| \big[W(t)\cdots W(\tau+1)W(\tau)\big]_{ij} - v_i(t) \big|\le c \mu^{t-\tau},
\end{align}
where $[\cdot]_{ij}$ denotes the $ij$th entry of a matrix. 
It has been shown in \cite{nedic} that $c=4$ and $\mu=(1-\frac{1}{n^{nL}})^{1/L}$.

To proceed, we define a time-dependent $n\times n$ matrix $S(t)$ whose $ij$th entry is 
\begin{equation}\label{eq:s}
    s_{ij}(t)\dfb \frac{w_{ij}(t)y_j(t)}{y_i(t+1)}=\frac{w_{ij}(t)y_j(t)}{\sum_{k=1}^n w_{ik}(t)y_k(t)}.
\end{equation}
It is worth emphasizing that $S(t)$ is independent of $x(t)$.
The following lemma guarantees that $S(t)$ is well defined.

\begin{lemma}\label{lemma:y_bound}
    If $\{ \bbb{G}(t) \}$ is uniformly strongly connected, then there exists a constant $\eta>0$ such that $n \ge y_i(t) \ge \eta$ for all $i$ and $t$.
\end{lemma}

It is easy to prove that $n\ge y_i(t)$. 
The lemma is essentially the same as Corollary~2~(b) in~\cite{nedic}, which further proves that if $\{ \bbb{G}(t) \}$ is uniformly strongly connected by sub-sequences of length $L$, then 
$\eta \ge \frac{1}{n^{nL}}$.




Define $z_i(t) \dfb x_i(t)/y_i(t)$ for each $i\in\mathcal{V}$. Then, 
\begin{align}
z_i(t+1) 
&= \frac{x_i(t+1)}{y_i(t+1)} 
= \frac{\sum_{j=1}^n w_{ij}(t)x_j(t)}{\sum_{j=1}^n w_{ij}(t)y_j(t)}
= \sum_{j=1}^n \frac{w_{ij}(t)x_j(t)}{\sum_{k=1}^n w_{ik}(t)y_k(t)}\nonumber\\
&= \sum_{j=1}^n \bigg[\frac{w_{ij}(t)y_j(t)}{\sum_{k=1}^n w_{ik}(t)y_k(t)}\bigg]z_j(t) 
= \sum_{j=1}^n s_{ij}(t)z_j(t),\label{eq:update_ratio}
\end{align}
which implies that $z(t+1)=S(t)z(t)$ where $z(t)=[z_1(t) \; z_2(t)\; \cdots \; z_n(t)]^\top \in\R^{n\times d}$ is the stack of all $z_i^\top(t)$, $i\in\scr V$. 
Actually $S(t)$ is always a stochastic matrix, as we will show shortly. 

Similar to the discrete-time state transition matrix, 
let $\Phi_W(t,\tau)\dfb W(t-1)\cdots W(\tau)$ with $t>\tau$,
and similarly, let $\Phi_S(t,\tau)\dfb S(t-1)\cdots S(\tau)$ with $t>\tau$.


\begin{lemma}\label{lemma:yixuan}
For $i,j\in\scr V$ and $t > \tau \ge 0$, there holds
$[\Phi_S(t,\tau)]_{ij}y_i(t)=[\Phi_W(t,\tau)]_{ij}y_j(\tau).$
\end{lemma}

\noindent
{\bf Proof of Lemma~\ref{lemma:yixuan}:}
The claim will be proved by induction on $t$. 
For the basis step, suppose that $t=\tau+1$. Then, from \eqref{eq:s},
\begin{align*}
[\Phi_S(\tau+1,\tau) ]_{ij} &= s_{ij}(\tau) = \frac{ y_j(\tau)  w_{ij}(\tau)}{y_i(\tau+1)} = \frac{ y_j(\tau)}{y_i(\tau+1)} [\Phi_W(\tau+1,\tau) ]_{ij}.
\end{align*}
Thus, in this case the claim is true.
For the inductive step, 
suppose that the claim holds for $t=h > \tau$, where $h$ is a positive integer, and that $t=h+1$. Then,
\begin{align*}
[\Phi_S(h+1,\tau)]_{ij} 
&= \sum_{k=1}^n  s_{ik}(h)  [\Phi_S(h,\tau)]_{kj} 
= \sum_{k=1}^n  \frac{w_{ik}(h)y_k(h)}{y_i(h+1)}   \frac{ y_j(\tau)  }{y_k(h)} [\Phi_W(h,\tau) ]_{kj}  \\
&= \frac{ y_j(\tau)  }{y_i(h+1)} \sum_{k=1}^n w_{ik}(h)   [\Phi_W(h,\tau) ]_{kj} 
= \frac{ y_j(\tau)  }{y_i(h+1)} [\Phi_W(h+1,\tau) ]_{ij} , 
\end{align*}
which establishes the claim by induction.
\hfill$\qed$

More can be said. 

\begin{lemma}\label{yixuan}
    If $\{ \bbb{G}(t) \}$ is uniformly strongly connected, then for any fixed $\tau\ge0$, $S(t)\cdots S(\tau+1)S(\tau)$ will converge to $\frac{1}{n} \1 y^\top(\tau)$ as $t\rightarrow\infty$. 
\end{lemma}





{\bf Proof of Lemma~\ref{yixuan}:}
From Lemma~\ref{lemma:pushsum_product}, for any given $\tau\ge 0$, there holds $\lim_{t\to\infty} [\Phi_W(t,\tau)] = v(\tau,\infty) \1^\top $, with the understanding that $v(\tau,\infty)$ is not necessarily a constant vector. From Lemma~\ref{lemma:yixuan} and the fact that $y(t) =  \Phi_W(t,\tau) y(\tau)$ for all $t > \tau$, 
\begin{align*}
\lim_{t\to\infty} [\Phi_S(t,\tau)]_{ij}  
&= \lim_{t\to\infty} \frac{ y_j(\tau)  }{y_i(t)} [\Phi_W(t,\tau) ]_{ij} 
= \lim_{t\to\infty} \frac{ y_j(\tau) [\Phi_W(t,\tau) ]_{ij} }{\sum_{k=1}^n [\Phi_W(t,\tau) ]_{ik} y_k(\tau)} \\
&=  \frac{ y_j(\tau) \lim_{t\to\infty} [\Phi_W(t,\tau) ]_{ij} }{\lim_{t\to\infty} \sum_{k=1}^n [\Phi_W(t,\tau) ]_{ik} y_k(\tau)} 
=  \frac{ y_j(\tau) v_i(\tau,\infty) }{  \sum_{k=1}^n v_i(\tau,\infty) y_k(\tau)}\\ 
&\overset{(a)}{=}  \frac{ y_j(\tau) }{  \sum_{k=1}^n y_k(\tau)} 
\overset{(b)}{=} \frac{y_j(\tau)}{n},
\end{align*}
where in (a) we used the fact that $v(\tau,\infty)>\0$ by Lemma~\ref{lemma:pushsum_product} and in (b) we used the fact that $\sum_{i=1}^n y_i(t) = n $ for all $t \ge 0$.
\hfill$\qed$

The limit in Lemma \ref{yixuan} is approached exponentially fast.

\begin{proposition} \label{lemma:expconvergen} 
    If $\{ \bbb{G}(t) \}$ is uniformly strongly connected, then for any fixed $\tau\ge0$, $S(t)\cdots S(\tau+1)S(\tau)$ will converge to $\frac{1}{n} \1 y^\top(\tau)$ exponentially fast as $t\rightarrow\infty$. 
\end{proposition}


{\bf Proof of Proposition~\ref{lemma:expconvergen}:}
From \eqref{mu}, there exist constants $c>0$ and $\mu \in [0,1)$ and a sequence of stochastic vectors $ \{ v(t)\}$ such that 
$
    | [\Phi_W(t+1,\tau)]_{ij} - v(t) |\le c \mu^{t-\tau}
$ for all $i,j \in \scr V$ and $t \ge \tau \ge 0$.
Recall that $\sum_{i=1}^n y_i(t)$ always equals $n$ and, by Lemma~\ref{lemma:y_bound},  all $y_i(t)$ are always positive. From Lemma~\ref{lemma:yixuan}, for all $t\ge\tau\ge0$,
\begin{align*}
&\Big| [\Phi_S(t+1,\tau)]_{ij} - \frac{y_j(\tau)}{n} \Big|
= \Big| \frac{ y_j(\tau) [\Phi_W(t+1,\tau) ]_{ij} }{y_i(t+1)}  - \frac{y_j(\tau)}{n} \Big|\\
=\;&\Big| \frac{ ny_j(\tau) [\Phi_W(t+1,\tau) ]_{ij} - y_j(\tau) [\Phi_W(t+1,\tau)y(\tau)]_i}{ny_i(t+1)}   \Big|
\\
=\;& \Big| \frac{n y_j(\tau) \left([\Phi_W(t+1,\tau)]_{ij} - v_i(t) + v_i(t) \right) - y_j(\tau)\sum_{k=1}^n ([\Phi_W(t+1,\tau)]_{ik} - v_i(t) + v_i(t) ) y_k(\tau) }{ny_i(t+1)} \Big|\\
=\;&  \Big| \frac{ n y_j(\tau) ([\Phi_W(t+1,\tau)]_{ij} - v_i(t) ) - y_j(\tau) \sum_{k=1}^n ( [\Phi_W(t+1,\tau)]_{ik}- v_i(t) ) y_k(\tau)    }{ ny_i(t+1) } \Big|\\
\le \;&   \frac{  n y_j(\tau) \left|[\Phi_W(t+1,\tau)]_{ij} - v_i(t) \right| +  y_j(\tau) \sum_{k=1}^n \left| [\Phi_W(t+1,\tau)]_{ik}- v_i(t) \right| y_k(\tau) }{ ny_i(t+1) } \\
\le \;&  \frac{ n y_j(\tau) c \mu^{t-\tau} + y_j(\tau) \sum_{k=1}^n c \mu^{t-\tau} y_k(\tau)     }{ ny_i(t+1) } 
= \frac{ 2 y_j(\tau) c \mu^{t-\tau} }{ y_i(t+1) } \le  \frac{2cn}{ \eta } \mu^{t-\tau},
\end{align*}
where we used Lemma~\ref{lemma:y_bound} in the last inequality. The above immediately implies the proposition.
\hfill$\qed$

The proposition immediately implies the following result and Theorem~\ref{thm:pushsum}.

\begin{corollary}\label{coro:sproduct}
    If $\{ \bbb{G}(t) \}$ is uniformly strongly connected, then $S(t)\cdots S(1)S(0)$ will converge to $\frac{1}{n}\1\1^\top$ exponentially fast as $t\rightarrow\infty$.\footnote{We were aware of this assertion for the first time from C\'esar A. Uribe in a private communication.}
\end{corollary}

{\bf Proof of Corollary \ref{coro:sproduct}:} The corollary is a special case of Proposition~\ref{lemma:expconvergen} by setting $\tau=0$. 
\hfill$\qed$

{\bf Proof of Theorem \ref{thm:pushsum}:}
From \eqref{eq:update_ratio}, $z(t+1)=S(t)z(t)=S(t)\cdots S(1)S(0)z(0)=S(t)\cdots S(1)S(0)x(0)$. From Corollary~\ref{coro:sproduct}, $z(t+1)$ will converge to $\frac{1}{n}\1\1^\top x(0) = (\frac{1}{n}\sum_{i=1}^nx_i(0))\1$ exponentially fast as $t\rightarrow\infty$, which completes the proof.
\hfill$\qed$

From the proof of Proposition~\ref{lemma:expconvergen}, each entry of $S(t)\cdots S(1)S(0)$ converges to $1/n$ exponentially fast at a rate of at least $\mu$, so does the convergence of $z_i(t)$ to $\frac{1}{n}\sum_{i=1}^nx_i(0)$. This convergence rate is consistent with that of the convergence of $W(t)\cdots W(1)W(0)$; see \eqref{mu}. Thus, our new proof and the conventional proof (see e.g. \cite{push,weighted,wallerton}) of Theorem~\ref{thm:pushsum} yield the same upper bound of the convergence rate of the push-sum algorithm.

Although the above proof of Theorem~\ref{thm:pushsum} looks more complicated than the conventional proof in the literature, it yields the following key property of the push-sum algorithm.

\subsection{A Key Property}\label{sec:key}

To proceed, we rewrite the push-sum algorithm in a different form which directly characterizes the dynamics of $z_i(t)=x_i(t)/y_i(t)$. 
From \eqref{eq:s} and \eqref{eq:update_ratio}, $z_i(t+1)=\sum_{j=1}^n s_{ij}(t)z_j(t)$ and $s_{ij}(t)$ satisfies the following assumption.

\begin{assumption}\label{assum:rowstochastic}
     There exists a constant $\gamma>0$ such that for all $i,j\in\scr V$ and $t$, $s_{ii}(t) \ge \gamma$ and $s_{ij}(t) \ge \gamma$ whenever $s_{ij}(t)>0$. For all $i\in\scr V$ and $t$, $\sum_{j=1}^n s_{ij}(t) = 1$.
\end{assumption}

\begin{lemma}\label{lemma:smatrices}
    If Assumption~\ref{assum:weighted matrix} holds, then $s_{ij}(t)$ satisfies Assumption~\ref{assum:rowstochastic} for each $t\ge 0$.
\end{lemma}

{\bf Proof of Lemma~\ref{lemma:smatrices}:}
From Assumption~\ref{assum:weighted matrix}, each $ W(t)$ is a column stochastic matrix whose diagonal entries are all positive and $ w_{ij}(t) \ge \beta$ whenever $ w_{ij}(t) > 0$. From \eqref{eq:s}, $s_{ij}(t)>0$ only if $w_{ij}(t)>0$. From Lemma~\ref{lemma:y_bound}, when $w_{ij}(t)>0$, 
$$ s_{ij}(t)=\frac{w_{ij}(t)y_j(t)}{y_i(t+1)}
\ge  \frac{ \beta \eta}{n}.$$
The above inequality and Assumption~\ref{assum:weighted matrix} imply that $s_{ij}(t)$ satisfies the first sentence of Assumption~\ref{assum:rowstochastic} with $\gamma=\beta\eta/n$.
For the second sentence of Assumption~\ref{assum:rowstochastic}, it is easy to see that 
\begin{align*}
    \sum_{j=1}^n s_{ij}(t) & = \sum_{j=1}^n \frac{w_{ij}(t)y_j(t)}{\sum_{k=1}^n w_{ik}(t)y_k(t)}  = 1
\end{align*}
for all $i\in\scr V$ and $t$, which completes the proof.
\hfill$\qed$

From Lemma~\ref{lemma:smatrices}, each $S(t)$ is a row stochastic matrix whose diagonal entries are all positive and whose nonzero entries are all uniformly bounded below by some positive number. More can be said. The following lemma shows that each $S(t)$ is compliant with the neighbor graph $\bbb{G}(t)$.

\begin{lemma}\label{lemma:sgraph}
    The graph of $S(t)$ is the same as the graph of $W(t)$ for all $t$.\footnote{The graph of an $n\times n$ matrix is a direct graph with $n$ vertices and an arc from vertex $i$ to vertex $j$ whenever the $ji$th entry of the matrix is~nonzero.}
\end{lemma}

{\bf Proof of Lemma~\ref{lemma:sgraph}:}
From \eqref{eq:s} and  Lemma~\ref{lemma:y_bound}, it is easy to see that $s_{ij}(t)>0$ if and only if $w_{ij}(t)$, which proves the lemma.
\hfill$\qed$

From~\eqref{eq:update_ratio}, $z(t+1)=S(t)z(t)$. 
The above lemmas imply that the dynamics of $z(t)$ is a nonlinear consensus process as $S(t)$ is dependent on $z(t)$. 
Such a transition in analysis from $x(t)$ dynamics to $z(t)$ dynamics has been used in \cite{iutzeler2013analysis}. 
To analyze such a process, we appeal to the following concept. To our knowledge, the concept has never been used to analyze the push-sum algorithm and push-sum based distributed algorithms.

\begin{definition}\label{def: absolute prob}
Let $\{ S(t) \}$ be a sequence of stochastic matrices. A sequence of stochastic vectors $\{ \pi(t) \}$ is an absolute probability sequence for $\{ S(t) \}$ if
$\pi^\top(t) = \pi^\top(t+1) S(t)$ for all $t\ge0$.
\end{definition}

This definition was first introduced by Kolmogorov who proved that every sequence of stochastic matrices has an absolute probability sequence \cite{kolmogorov}. An alternative proof of this fact was given by Blackwell \cite{blackwell}.
In general, a sequence of stochastic matrices may have more than one absolute probability sequence; when the sequence of stochastic matrices is ``ergodic'',\footnote{
A sequence of stochastic matrices $\{S(t)\}$ is called ergodic if $\lim_{t\rightarrow\infty}S(t)\cdots S(\tau+1)S(\tau)=\1v^\top(\tau)$ for all $\tau$, where each $v(\tau)$ is a stochastic vector.} it has a unique absolute probability sequence \cite[Lemma~1]{tacrate}. It is easy to see that when $S(t)$ is a fixed irreducible stochastic matrix $S$, $\pi(t)$ is simply the normalized left eigenvector of $S$ for eigenvalue one, and when $\{S(t)\}$ is an ergodic sequence of doubly stochastic matrices, $\pi(t)=(1/n)\1$.
More can be said.

\begin{lemma} \label{lemma:bound_pi_jointly}
    {\rm (Theorem 4.8 in \cite{touri2012product})}
    Let $\{S(t)\}$ be a sequence of stochastic matrices satisfying Assumption~\ref{assum:rowstochastic}. If the graph sequence of $\{\bbb{G}(t)\}$ is uniformly strongly connected, then there exists a unique absolute probability sequence $\{ \pi(t) \}$ for the matrix sequence $\{S(t)\}$ and a constant $\pi_{\min} \in (0,1)$ such that $\pi_i(t) \ge \pi_{\min}$ for all $i$ and $t$.
\end{lemma}

A particular important property of the absolute probability sequence for $\{S(t)\}$ is as follows.

\begin{proposition} \label{lemma:push-sum_pi_intfty}
    If $\{ \bbb{G}(t) \}$ is uniformly strongly connected, then the sequence of stochastic matrices $\{S(t)\}$ has a unique absolute probability sequence $\{\pi(t)\}$ with
    $\pi_i(t)=\frac{y_i(t)}{n}$ for all $i\in\scr V$ and $t \ge 0$. 
\end{proposition}


The proposition is a consequence of Lemma 1 in \cite{tacrate}. We provide two alternative proofs. 

{\bf Proof of Proposition~\ref{lemma:push-sum_pi_intfty}:}
First, Lemma~\ref{yixuan} shows that $\{S(t)\}$ is ergodic, so it must have a unique absolute probability sequence $\{\pi(t)\}$. 
From Definition~\ref{def: absolute prob} and Lemma~\ref{yixuan}, for any $\tau \ge 0$,
\begin{align*}
\pi^\top(\tau) 
&= \pi^\top(\tau+1) S(\tau) 
= \pi^\top(\tau+2) S(\tau+1)S(\tau)
= \cdots = \lim_{t\to\infty}\pi^\top(t+1) S(t)\cdots S(\tau+1)S(\tau) \\
&= \lim_{t\to\infty} \frac{1}{n}\pi^\top(t+1)  \1 y^\top(\tau) 
= \frac{1}{n} y^\top(\tau),
\end{align*}
which proves the statement. 

Alternatively, we can also prove the proposition by showing that the sequence $\{\pi(t)\}$ with 
$\pi_i(t) = \frac{y_i(t)}{n}$ satisfies $\pi^\top(t) = \pi^\top(t+1) S(t)$. To see this, from \eqref{eq:s} and Assumption~\ref{assum:weighted matrix}, 
\begin{align*}
    [\pi^\top(t+1) S(t)]_j 
    = \sum_{i=1}^n \frac{y_i(t+1)}{n} s_{ij}(t)  
    = \sum_{i=1}^n \frac{y_i(t+1)}{n} \frac{w_{ij}(t)y_j(t)}{y_i(t+1)} 
    = \sum_{i=1}^n  \frac{w_{ij}(t)y_j(t)}{n} 
    = \frac{y_j(t)}{n} = \pi_j(t) 
\end{align*}
for all $j \in \scr V$, 
where $[\cdot]_{j}$ denotes the $j$th entry of a vector. 
\hfill$\qed$

Proposition~\ref{lemma:push-sum_pi_intfty} immediately implies that 
$y^\top(t) = y^\top(t+1) S(t)$ for all $t$. 
Proposition~\ref{lemma:push-sum_pi_intfty} and Lemma~\ref{lemma:y_bound} imply that $\{\pi(t)\}$ of $\{S(t)\}$ is always positive and its entries are bounded below by $\pi_{\min}=\eta/n$.

In next sections, we will illustrate the usefulness of the above key property by applying it to push-sum based distributed optimization. Specifically, we appeal to the property to construct a novel time-varying Lyapunov function for distributed convex optimization which yields improved and optimal convergence rates of the so-called subgradient-push and stochastic gradient-push algorithms.

Since the stochastic matrix sequence $S(t)$ defined by \eqref{eq:s} is purely based on the $y_i(t)$ variables and is thus
independent of the $x_i(t)$ variables of the push-sum algorithm, so its absolute probability sequence. Considering the fact that the push-sum algorithm \eqref{pushsumx}--\eqref{pushsumy}, subgradient-push algorithm \eqref{eq:pushsub_x}--\eqref{eq:pushsub_y}, and stochastic gradient-push \eqref{eq:spushsub_x}--\eqref{eq:spushsub_y} share the same $y_i(t)$ dynamics which is independent of their $x_i(t)$ dynamics, all the results of $\{S(t)\}$ and its absolute probability sequence derived in this subsection apply to the subgradient-push and stochastic gradient-push algorithms. 
This is also the case for any other push-sum based distributed algorithms with independent $y_i(t)$ dynamics. 

\section{Subgradient-Push}\label{sec:subgrad}

In this section, we revisit the subgradient-push algorithm first proposed in \cite{nedic} and analyze its convergence rate using the key property of push-sum proved in Section~\ref{sec:key}. Our novel analysis tool will improve the convergence rate of the subgradient-push algorithm from $O(\ln t/\sqrt{t})$ to $O(1/\sqrt{t})$.

Consider the same network of $n$ agents as described at the beginning of Section~\ref{sec:pushsum}. The goal of the $n$ agents is to cooperatively to minimize the cost function
$$f(z)=\frac{1}{n}\sum_{i=1}^n f_i(z),$$ where each $f_i$ is a ``private'' convex (not necessarily differentiable) cost function only known to agent $i$. It is assumed that the set of optimal solutions to $f$, denoted by $\scr Z$, is nonempty.

Since each $f_i$ is not necessarily differentiable, the gradient descent method may not be applicable. Instead, the subgradient method \cite{subgradient} can be applied. For a convex function $h : \R^d\rightarrow \R$, a vector $g\in\R^d$ is called a subgradient of $h$ at point $x$ if
\eq{
h(y)\ge h(x) + g^\top (y-x) \;\; {\rm for \; all} \;\; y\in\R^d.
\label{eq:subgradient}}
Such a vector $g$ always exists and may not be unique. In the case when $h$ is differentiable at point $x$, the subgradient $g$ is unique and equals $\nabla h(x)$, the gradient of $h$ at $x$. Thus, the subgradient can be viewed as a generalization of
the notion of the gradient. From \eqref{eq:subgradient} and the Cauchy-Schwarz inequality, for any $x,y\in\R^d$,
\eq{
h(y) - h(x) \ge - G \| y-x\|,
\label{eq:G}}
where $G$ is an upper bound for the 2-norm of the subgradients of $h$ at both $x$ and $y$. We use $\|\cdot\|$ to denote the (induced) 2-norm throughout the paper.

The subgradient method was first proposed in \cite{subgradient} and the first distributed subgraident method was proposed in \cite{nedic2009distributed}, which is based on average consensus (or distributed averaging).
The subgradient-push algorithm, proposed in \cite{nedic},  is as follows\footnote{The algorithm is written in a different but mathematically equivalent form in \cite{nedic}.}:
\begin{align}
    x_i(t+1) &= \sum_{j\in\scr{N}_i(t)} w_{ij}(t)\Big[x_j(t) - \alpha(t)g_j(t)\Big], \;\;\;\;\; x_i(0)\in\R^d,\label{eq:pushsub_x}\\  y_i(t+1) &= \sum_{j\in\scr{N}_i(t)} w_{ij}(t)y_j(t),\;\;\;\;\; y_i(0)=1\label{eq:pushsub_y},
\end{align}
where $\alpha(t)$ is the stepsize, $g_j(t)$ is a subgradient of $f_j(z)$ at $x_j(t)/y_j(t)$, and $w_{ij}(t)$, $j\in\scr{N}_i(t)$, are positive weights satisfying Assumption \ref{assum:weighted matrix}.  

In implementation, at each time $t$, each agent $j$ transmits two pieces of information, $w_{ij}(t)[x_j(t) - \alpha(t)g_j(t)]$ and $w_{ij}(t)y_j(t)$, to its out-neighbour $i$, and then each agent $i$ updates its two variables as above. Note that if all $\alpha(t)g_j(t) = 0$, the algorithm simplifies to the push-sum algorithm \eqref{pushsumx}--\eqref{pushsumy}. Thus, at each time, each agent first performs a subgradient operation, and then follows the push-sum updates. This is why the algorithm \eqref{eq:pushsub_x}--\eqref{eq:pushsub_y} is called subgradient-push. In Section~\ref{sec:mix} we will study ``push-subgradient'' in which the order of subgradient and push-sum operations is swapped.

To state the convergence results of the subgradient-push algorithm \eqref{eq:pushsub_x}--\eqref{eq:pushsub_y}, we need the following assumption and notation.

\begin{assumption} \label{assum:step-size}
    The step-size sequence $\{\alpha(t)\}$ is positive, non-increasing, and satisfies $\sum_{t=0}^\infty \alpha(t) = \infty$ and $\sum_{t=0}^\infty \alpha^2(t) < \infty$. 
\end{assumption}

Recall $z_i(t)=x_i(t)/y_i(t)$ and define $\bar z(t) \dfb \frac{1}{n}\sum_{i=1}^n z_i(t)$.

\begin{theorem} \label{thm:bound_everage_n_convex}
    Suppose that $\{ \bbb{G}(t) \}$ is uniformly strongly connected and that $\|g_i(t)\|$ is uniformly bounded for all $i$ and~$t$.
\begin{itemize}
    \item[1)] If the stepsize $\alpha(t)$ is time-varying and satisfies Assumption~\ref{assum:step-size}, then with $z^* \in \scr{Z}$, 
\begin{align*}
    \lim_{t\rightarrow\infty}f\bigg(\frac{\sum_{\tau =0}^t \alpha(\tau) \bar z(\tau) }{\sum_{\tau =0}^t \alpha(\tau)}\bigg) = f(z^*).
\end{align*}
    \item[2)] If the stepsize is fixed and  $\alpha(t) = 1/\sqrt{T}$ for $T>0$ steps, i.e., $t\in\{0,1,\ldots,T-1\}$, then 
\begin{align*}
    f\bigg(\frac{\sum_{\tau =0}^{T-1} \bar z(\tau)  }{ T } \bigg) - f(z^*) 
    & \le O\Big(\frac{1}{ \sqrt{T}}\Big).
\end{align*}
\end{itemize}
\end{theorem}

The above theorem establishes the convergence rate of $f((\sum_{\tau =0}^{T-1} \bar z(\tau))/T)$, as conventionally did in average consensus based subgradient \cite{nedic2009distributed}, and the rate is of $O(1/\sqrt{t})$, which is the same as that of the conventional single-agent subgradient method \cite[Theorem 7]{nedic2018network}. Thus, the derived convergence rate is optimal.

Theorem \ref{thm:bound_everage_n_convex} is actually a consequence of the following refined result, which further provides finite-time error bounds for the subgradient-push algorithm \eqref{eq:pushsub_x}--\eqref{eq:pushsub_y}.

\begin{theorem} \label{thm:bound_everage_n_convex_bound}
    Suppose that $\{ \bbb{G}(t) \}$ is uniformly strongly connected by sub-sequences of length $L$ and that $\|g_i(t)\|$ is uniformly bounded above by a positive number $G$ for all $i$ and $t$.
\begin{itemize}
    \item[1)] 
If the stepsize $\alpha(t)$ is time-varying and satisfies Assumption~\ref{assum:step-size}, then for all $t \ge 0$,
\begin{align}
    f\bigg(\frac{\sum_{\tau =0}^t \alpha(\tau) \bar z(\tau) }{\sum_{\tau =0}^t \alpha(\tau)}\bigg) - f(z^*) 
     \le \;& \frac{ \| \bar z(0) - z^* \|^2 + G^2\sum_{\tau =0}^t \alpha^2(\tau)}{2\sum_{\tau =0}^t \alpha(\tau) }  
    + \frac{2G\alpha(0) \sum_{i=1}^n \|\bar z(0) - z_i(0) \| }{n\sum_{\tau =0}^t \alpha(\tau) } \nonumber\\
    & + \frac{32G}{\eta} \cdot \frac{\sum_{\tau =0}^{t-1} \alpha(\tau) \mu^\tau  }{\sum_{\tau =0}^t \alpha(\tau) } \sum_{i=1}^n \big\| x_i(0) - \alpha(0) g_i(0) \big\|  \nonumber\\
    & + \frac{32nG^2}{\eta(1-\mu)} \cdot\frac{\sum_{\tau =0}^{t-1} \alpha(\tau)   \big(  \alpha(0) \mu^{\frac{\tau}{2}}  +  \alpha(\ceil{\frac{\tau}{2}}) \big)}{\sum_{\tau =0}^t \alpha(\tau) }\label{eq:bound_timevarying}.
\end{align} 
\item[2)] If the stepsize is fixed and  $\alpha(t) = 1/\sqrt{T}$ for $T>0$ steps, then
\begin{align}
    f\bigg(\frac{\sum_{\tau =0}^{T-1} \bar z(\tau)  }{ T } \bigg) - f(z^*) 
     \le \;& \frac{2G \sum_{i=1}^n \|\bar z(0) - z_i(0) \|}{nT } +\frac{ \| \bar z(0) - z^* \|^2 + G^2 }{ 2\sqrt{T} } \nonumber\\
    & + \frac{32 G }{\eta(1 - \mu)T}  \sum_{i=1}^n \big\| x_i(0) - \frac{1}{\sqrt{T}} g_i(0) \big\|  + \frac{32n G^2 }{ \eta(1- \mu)\sqrt{T}} \label{eq:bound_fixed}.
\end{align}
\end{itemize}
Here $\eta$ and $\mu$ are positive constants which satisfy $\eta \ge \frac{1}{n^{nL}}$ and $\mu \le (1-\frac{1}{n^{nL}})^{1/L}$, respectively, and $\ceil{\cdot}$ denotes the ceiling function.
\end{theorem}


The above theorem characterizes convergence rates for a network-wide averaged state. The following theorem provides convergence rates for each individual agent in the subgradient-push algorithm.

\begin{theorem} \label{thm:bound_everage_zi}
    Suppose that $\{ \bbb{G}(t) \}$ is uniformly strongly connected by sub-sequences of length $L$ and that $\|g_i(t)\|$ is uniformly bounded above by a positive number $G$ for all $i$ and $t$.
\begin{itemize}
    \item[1)] 
If the stepsize $\alpha(t)$ is time-varying and satisfies Assumption~\ref{assum:step-size}, then for all $t \ge 0$ and $k\in\scr{V}$,
\begin{align*}
    f\bigg(\frac{\sum_{\tau =0}^t \alpha(\tau) z_k(\tau) }{\sum_{\tau =0}^t \alpha(\tau)}\bigg) - f(z^*) 
    \le \; & \frac{ \| \bar z(0) - z^* \|^2 + G^2\sum_{\tau =0}^t \alpha^2(\tau)}{2\sum_{\tau =0}^t \alpha(\tau) }  \nonumber\\
    &  + \frac{G\alpha(0) \sum_{i=1}^n \big(\|\bar z(0) - z_i(0) \| + \|z_k(0)  - z_i(0) \|\big)}{n\sum_{\tau =0}^t \alpha(\tau) } \nonumber\\
     & + \frac{32G}{\eta} \cdot\frac{\sum_{\tau =0}^{t-1} \alpha(\tau) \mu^\tau  }{\sum_{\tau =0}^t \alpha(\tau) } \sum_{i=1}^n \big\| x_i(0) - \alpha(0) g_i(0) \big\|  \nonumber\\
     & + \frac{32nG^2}{\eta(1-\mu)} \cdot \frac{\sum_{\tau =0}^{t-1} \alpha(\tau)   \big(  \alpha(0) \mu^{\frac{\tau}{2}}  +  \alpha(\ceil{\frac{\tau}{2}}) \big)}{\sum_{\tau =0}^t \alpha(\tau) }. 
\end{align*} 
\item[2)] If the stepsize is fixed and  $\alpha(t) = 1/\sqrt{T}$ for $T>0$ steps, then for any $k\in\scr{V}$,
\begin{align*}
    f\bigg(\frac{\sum_{\tau =0}^{T-1} z_k(\tau)  }{ T } \bigg) - f(z^*) 
    \le \; & \frac{ \| \bar z(0) - z^* \|^2 + G^2 }{ 2\sqrt{T} }+ \frac{32n G^2 }{ \eta(1- \mu)\sqrt{T}} \nonumber\\
    & + \frac{G \sum_{i=1}^n \big(\|\bar z(0)  - z_i(0) \|+\|z_k(0) - z_i(0) \|\big)}{nT } \nonumber\\
    & + \frac{32 G }{\eta(1 - \mu)T}  \sum_{i=1}^n \big\| x_i(0) - \frac{1}{\sqrt{T}} g_i(0) \big\|.  
\end{align*}
\end{itemize}
Here the positive constants $\eta$ and $\mu<1$ are the same as in Theorem \ref{thm:bound_everage_n_convex_bound}.
\end{theorem}

Using the same argument as in the proof
of Theorem \ref{thm:bound_everage_n_convex}, we have for each agent $k\in\scr V$, with a time-varying stepsize $\alpha(t)$ satisfying Assumption~\ref{assum:step-size},
\begin{align*}
    \lim_{t\rightarrow\infty}f\bigg(\frac{\sum_{\tau =0}^t \alpha(\tau) z_k(\tau) }{\sum_{\tau =0}^t \alpha(\tau)}\bigg) = f(z^*),
\end{align*}
and with a fixed stepsize $\alpha(t) = 1/\sqrt{T}$ for $T>0$ steps, 
\begin{align*}
    f\bigg(\frac{\sum_{\tau =0}^{T-1} z_k(\tau)  }{ T } \bigg) - f(z^*) 
    & \le O\Big(\frac{1}{ \sqrt{T}}\Big).
    \end{align*}

\subsection{Analysis}

In this subsection, we provide a novel analysis of the subgradient-push algorithm \eqref{eq:pushsub_x}--\eqref{eq:pushsub_y} and proofs of Theorems \ref{thm:bound_everage_n_convex}--\ref{thm:bound_everage_zi}.
Before moving to the long proofs, we first highlight the technical challenge and elaborate our corresponding technical contribution. 

The key challenge in analyzing (sub)gradient-push is to find a suitable Lyapunov function. For distributed (sub)gradient using the averaging scheme with doubly stochastic matrices, there is no $y$ state dynamics and a common Lyapunov function is $\|\frac{1}{n}\sum_{i=1}^n x_i(t)-x^*\|^2$ (see e.g. \cite{nedic2009distributed,yuan2016convergence}). A natural guess for (sub)gradient-push would be $\|\frac{1}{n}\sum_{i=1}^n z_i(t)-z^*\|^2$, but it does not work (so far). The work of \cite{nedic} makes use of an untypical one, based on a complicated convex combination of agents' states, yet it yields a non-optimal convergence rate, so does \cite{nedic2016stochastic}.  
In the sequel, we will construct a novel and simple Lyapunov function tailored for push-sum based algorithms.

We first rewrite the subgradient-push algorithm as follows. 
From \eqref{eq:pushsub_x}--\eqref{eq:pushsub_y}, 
\begin{align}
z_i(t+1) 
&= \frac{x_i(t+1)}{y_i(t+1)} 
= \frac{\sum_{j=1}^n w_{ij}(t)[x_j(t)- \alpha(t)g_j(t)]}{\sum_{j=1}^n w_{ij}(t)y_j(t)}
= \sum_{j=1}^n \frac{w_{ij}(t)[x_j(t)- \alpha(t)g_j(t)]}{\sum_{k=1}^n w_{ik}(t)y_k(t)} \nonumber\\
&= \sum_{j=1}^n \bigg[\frac{w_{ij}(t)y_j(t)}{\sum_{k=1}^n w_{ik}(t)y_k(t)}\bigg]\bigg[z_j(t)- \alpha(t)
\frac{g_j(t)}{y_j(t)}\bigg]
= \sum_{j=1}^n s_{ij}(t)\bigg[z_j(t)- \alpha(t)
\frac{g_j(t)}{y_j(t)}\bigg],\label{eq:z_i(t+1)}
\end{align}
where $s_{ij}(t)$ is defined in \eqref{eq:s}.
In addition, 
\begin{align*}
    \bar z(t+1) 
    &= \frac{1}{n}\sum_{i=1}^n z_i(t+1) 
    = \frac{1}{n}\sum_{i=1}^n \sum_{j=1}^n s_{ij}(t)\bigg[z_j(t)- \alpha(t) \frac{g_j(t)}{y_j(t)}\bigg].
\end{align*}
To proceed, define the following time-dependent quantity: 
\eq{\langle z(t) \rangle \dfb \pi^\top(t)z(t) = \frac{1}{n} \sum_{i=1}^n y_i(t) z_i(t),\label{eq:<>}} 
where $\{\pi(t)\}$ is the unique absolute probability sequence of $\{S(t)\}$ and we used Proposition \ref{lemma:push-sum_pi_intfty} in the second equality. It is easy to see that the above quantity is a time-varying convex combination of all $z_i(t)$.   
Then, from \eqref{eq:z_i(t+1)} and Definition~\ref{def: absolute prob},
\begin{align}
    \langle z(t+1) \rangle 
    = \sum_{i=1}^n \pi_i({t+1}) z_i(t+1) 
    &= \sum_{i=1}^n \sum_{j=1}^n \pi_i({t+1})  s_{ij}(t)\bigg[z_j(t)- \alpha(t) \frac{g_j(t)}{y_j(t)}\bigg] \nonumber\\
    &= \sum_{j=1}^n \pi_j(t) \bigg[z_j(t)- \alpha(t) \frac{g_j(t)}{y_j(t)}\bigg]=\langle z(t) \rangle - \frac{\alpha(t)}{n}\sum_{i=1}^n g_i(t),\label{eq:<z>update}
\end{align}
where we use the Proposition~\ref{lemma:push-sum_pi_intfty} in the last equality.

The above iterative dynamics of $\langle z\rangle$ can be treated (though not exactly the same) as a single-agent subgradient process for the convex cost function $\frac{1}{n}\sum_{i=1}^n f_i(z)$. This is a critical intermediate step.

The remaining analysis logic is as follows. Using the inequality $\| z_i(t) - z^*\|^2 \le 2 \|\langle z(t) \rangle - z^*\|^2 + 2 \|\langle z(t) \rangle - z_i(t)\|^2$, the analysis is then to bound $\|\langle z(t) \rangle - z^*\|^2$ and $\|\langle z(t) \rangle - z_i(t)\|^2$ separately. For the term $\|\langle z(t) \rangle - z_i(t)\|^2$, since all $z_i$ form a consensus process and $\langle z(t) \rangle$ is always a convex combination of all $z_i(t)$, the term can be bounded using consensus related techniques and relatively easy to deal with. Our most analysis will focus on bounding the term $\|\langle z(t) \rangle - z^*\|^2$. It is worth noting that from \eqref{eq:<>},
$\|\langle z(t) \rangle - z^*\|^2=\|\frac{1}{n}\sum_{i=1}^n y_i(t)(z_i(t)-z^*)\|^2$,
which is a time-varying quadratic Lyapunov comparison function. This form helps us derive an optimal rate. 
The above Lyapunov function can be further simplified. Recall $z_i(t)=x_i(t)/y_i(t)$. It implies that the actual Lyapunov function for subgradient-push is $\|\frac{1}{n}\sum_{i=1}^n x_i(t)-z^*\|^2$, which is quite counterintuitive.
Note that this also implies that update \eqref{eq:<z>update} is equivalent to $\bar x(t+1)=\bar x(t) - \frac{\alpha(t)}{n}\sum_{i=1}^n g_i(t)$ where $\bar x(t) = \frac{1}{n}\sum_{i=1}^n x_i(t)$, which is almost the same as the case of average consensus based subgradient \cite{nedic2009distributed} except that each subgradient $g_i$ is taken at point $z_i$ instead of $x_i$. But this $\bar x$ dynamics is elusive without Proposition~\ref{lemma:push-sum_pi_intfty}.

To prove Theorems~\ref{thm:bound_everage_n_convex} and~\ref{thm:bound_everage_n_convex_bound}, we need the following lemma.

 \begin{lemma}\footnote{The authors thank Wenhan Gao for pointing out a mistake in the original version
of this lemma. The minor mistake appeared in the conference version \cite{pushcdc} which thus had incorrect but easily fixable upper bound expressions.} \label{lemma:bound_consensus_push_SA} 
    If $\{ \bbb{G}(t) \}$ is uniformly strongly connected by sub-sequences of length $L$ and $\|g_i(t)\|$ is uniformly bounded above by a positive number $G$ for all $i$ and $t$,
    then for all $t \ge 0$ and $i \in \mathcal{V}$,
 \begin{align*}
    \Big\| z_i(t+1) - \frac{1}{n} \sum_{k=1}^n \big(x_k(t) -\alpha(t) g_k(t)\big)\Big\| 
    \le \frac{8}{\eta}  \mu^t  \sum_{k=1}^n \| x_k(0) - \alpha(0) g_k(0) \| +  \frac{8nG}{\eta}  \sum_{s=1}^t \mu^{t-s} \alpha(s).
\end{align*}
If, in addition, Assumption~\ref{assum:step-size} holds, then for all $t \ge 0$ and $i \in \mathcal{V}$,
\begin{align*}
    \Big\| z_i(t+1) - \frac{1}{n} \sum_{k=1}^n \big(x_k(t) -\alpha(t) g_k(t)\big)\Big\| 
     \le \frac{8}{\eta}  \mu^t  \sum_{k=1}^n \| x_k(0) - \alpha(0) g_k(0) \| 
     +  \frac{8nG}{\eta(1-\mu)}   \big(  \alpha(0) \mu^{t/2}  +  \alpha(\ceil{t/2}) \big).
\end{align*}
Here $\eta>0$ and $\mu\in(0,1)$ are defined in Lemma~\ref{lemma:y_bound} and \eqref{mu}, respectively.
\end{lemma}

Lemma \ref{lemma:bound_consensus_push_SA} can be regarded as a high-dimensional extension of Lemma 1 and
Corollary 3 in \cite{nedic}, and its proof follows the ideas in the proofs of Lemma 1 and
Corollary 3 in \cite{nedic}.

To prove the lemma, we need the following notation. Let 
\begin{align*}
 g(t)\dfb \begin{bmatrix}g_1^\top(t)\cr\vdots\cr g_{n}^\top(t)\end{bmatrix}\in\R^{n\times d}.
\end{align*}
From \eqref{eq:pushsub_x}, it is easy to see that $x(t+1)=W(t)(x(t)-\alpha(t)g(t))$. We will use this fact without special mention in the sequel.

{\bf Proof of Lemma~\ref{lemma:bound_consensus_push_SA}:}
Let $h(t) = x(t) -\alpha(t) g(t) $. 
Then, 
\begin{align*}
    h(t+1) 
    &= x(t+1) -\alpha(t+1) g(t+1) 
    = W(t) h(t) -\alpha(t+1) g(t+1) \\
    & = \Phi_W(t+1,0) h(0) - \sum_{l=1}^t \alpha(l) \Phi_W(t+1,l) g(l) -\alpha(t+1) g(t+1).
\end{align*}
Note that
\begin{align}
    W(t+1) h(t+1)
     &= \Phi_W(t+2,0) h(0) - \sum_{l=1}^{t+1} \alpha(l) \Phi_W(t+2,l) g(l),\label{eq:h_update_tplus1}\\
    \1^\top h(t+1) 
    & = \1^\top h(0) - \sum_{l=1}^{t+1} \alpha(l) \1^\top g(l). \label{eq:h_update_1}
\end{align}
From Lemma~\ref{lemma:pushsum_product} and \eqref{mu}, there exists a sequence of stochastic vectors $ \{ \phi(t)\}$  such that for all $i,j \in \mathcal{V}$ and $t\ge s \ge 0$,
$| [\Phi_W(t+1,s)]_{ij} - \phi_i(t) |\le 4 \mu^{t-s}$.
Let 
\eq{D(s:t) \dfb \Phi_W(t+1,s) - \phi(t) \1^\top.\label{eq:D(s,t)}} 
From \eqref{eq:h_update_tplus1} and \eqref{eq:h_update_1}, 
\begin{align*}
    & W(t+1) h(t+1)  - \phi(t+1)\1^\top h(t+1)\\
    = \; & \Phi_W(t+2,0) h(0) - \sum_{l=1}^{t+1} \alpha(l) \Phi_W(t+2,l) g(l)
    - \phi(t+1) \Big(\1^\top h(0) - \sum_{l=1}^{t+1}  \alpha(l) \1^\top g(l)\Big) \\
    = \; & \big(\Phi_W(t+2,0) - \phi(t+1) \1^\top\big) h(0) 
    - \sum_{l=1}^{t+1} \alpha(l)\big( \Phi_W(t+2,l) - \phi(t+1) \1^\top\big) g(l) \\
    = \; & D(0:t+1) h(0) - \sum_{l=1}^{t+1} \alpha(l) D(l:t+1) g(l),
\end{align*}
which implies that
$x(t+1) = W(t) h(t) 
    = \phi(t)\1^\top h(t) + D(0:t) h(0) - \sum_{l=1}^{t} \alpha(l) D(l:t) g(l)$.
From \eqref{eq:pushsub_y} and \eqref{eq:D(s,t)},
$
    y(t+1) = \Phi_W(t+1,0) y(0) = D(0:t) \1 + n \phi(t)
$, or equivalently, $y_i(t+1) = [\Phi_W(t+1,0)\1]_i = [D(0:t) \1]_i + n \phi_i(t)$.
Thus, for all $i\in\scr V$,
\begin{align*}
    & z_i(t+1) - \frac{ h(t)^\top \1}{n} = \frac{x_i(t+1)}{y_i(t+1)} - \frac{ h(t)^\top \1}{n}\\
    = \; & \frac{\phi_i(t) h(t)^\top \1 + \sum_{k=1}^n [D(0:t)]_{ik} h_k(0) - \sum_{l=1}^{t} \alpha(l) \sum_{k=1}^n [D(l:t+1)]_{ik} g_k(l)}{[D(0:t) \1]_i + n \phi_i(t)} - \frac{ h(t)^\top \1}{n} \\
    = \; & \frac{n \sum_{k=1}^n [D(0:t)]_{ik} h_k(0)   - [D(0:t) \1]_i h(t)^\top \1 -  n\sum_{l=1}^{t} \alpha(l) \sum_{k=1}^n [D(l:t+1)]_{ik} g_k(l)   }{n[D(0:t) \1]_i + n^2 \phi_i(t)},
\end{align*}
where we use $h_i(t)\in\R^d$ to denote the $i$th column of $h^\top(t)\in\R^{d\times n}$.
From Lemma \ref{lemma:y_bound}, $y_i(t+1) \ge \eta$, so is $[D(0:t) \1]_i + n \phi_i(t)$. Thus, 
\begin{align}
    \Big\| z_i(t+1) - \frac{ h(t)^\top \1}{n} \Big\| 
    \le \; & \frac{n \|\sum_{k=1}^n [D(0:t)]_{ik} h_k(0)\|  + \|[D(0:t) \1]_i h(t)^\top \1 \| }{n[D(0:t) \1]_i + n^2 \phi_i(t)} \nonumber\\
    & + \frac{n\sum_{l=1}^{t} \alpha(l) \|\sum_{k=1}^n [D(l:t)]_{ik} g_k(l)\|  }{n[D(0:t) \1]_i + n^2 \phi_i(t)} \nonumber\\
    \le \; & \frac{n  (\max_k [D(0:t)]_{ik}) \sum_{k=1}^n  \| h_k(0)\|  + \|[D(0:t) \1]_i h(t)^\top \1 \| }{n[D(0:t) \1]_i + n^2 \phi_i(t)} \nonumber\\
    &  + \frac{n\sum_{l=1}^{t} \alpha(l) (\max_k [D(l:t)]_{ik}) \sum_{k=1}^n \| g_k(l)\|  }{n[D(0:t) \1]_i + n^2 \phi_i(t)} \nonumber\\
    \le \; & \frac{1}{n \eta} \Big[ n \big(\max_j | [D(0:t)]_{ij} |\big)  \sum_{k=1}^n \| h_k(0)\| \nonumber\\
    & + n\sum_{l=1}^{t} \alpha(l) \big(\max_j | [D(l:t)]_{ij} |\big)   \sum_{k=1}^n \| g_k(l)\|  +n \big(\max_j | [D(0:t)]_{ij} |\big)\| h(t)^\top \1 \|\Big] \nonumber\\
    \le \; & \frac{1}{ \eta} \Big[  4 \mu^{t} \sum_{k=1}^n  \|h_k(0)\| + \sum_{l=1}^{t} \alpha(l) 4 \mu^{t-l} \sum_{k=1}^n  \|g_k(l)\| 
    + 4 \mu^{t} \| h(t)^\top \1 \|\Big].\label{eq:middle}
\end{align}
Note that, from \eqref{eq:h_update_1}, 
$
    \| \1^\top h(t+1) \|
     \le \| \1^\top h(0)\| + \| \sum_{l=1}^{t+1} \alpha(l) \1^\top g(l)\| .
$
Then, from \eqref{eq:middle},
\begin{align*}
    & \Big\| z_i(t+1) - \frac{ h(t)^\top \1}{n} \Big\| \\
    \le \; & 
    \frac{4}{ \eta} \Big[   \mu^{t}  \sum_{k=1}^n \| h_k(0)\| + \sum_{l=1}^{t} \alpha(l)  \mu^{t-l} \sum_{k=1}^n  \|g_k(l)\|  +  \mu^{t} \| \1^\top h(0)\| +  \mu^{t} \Big\| \sum_{l=1}^{t} \alpha(l) \1^\top g(l)\Big\| \Big]\\
    \le \; & \frac{4}{ \eta} \Big[   2 \mu^{t}  \sum_{k=1}^n \| h_k(0)\| + 2 \sum_{l=0}^{t} \alpha(l)  \mu^{t-l} \sum_{k=1}^n \| g_k(l)\|  \Big] \\
    \le \; & \frac{8}{\eta}  \mu^t  \sum_{k=1}^n \| x_k(0) - \alpha(0) g_k(0) \| + \frac{8nG}{\eta} \sum_{l=0}^t  \mu^{t-l} \alpha(l).
\end{align*}
If the stepsize sequence $\{ \alpha(t) \}$ satisfies Assumption~\ref{assum:step-size}, the above inequality implies that
\begin{align*}
    \Big\| z_i(t+1) - \frac{ h(t)^\top \1}{n} \Big\| 
    & \le \frac{8}{\eta}  \mu^t  \sum_{k=1}^n \| x_k(0) - \alpha(0) g_k(0) \|  +  \frac{8nG}{\eta} \bigg(\sum_{l=0}^{\floor{t/2}}  \mu^{t-l} \alpha(l)
    + \sum_{l=\ceil{t/2}}^{t}  \mu^{t-l} \alpha(l) \bigg)   \\
    & \le \frac{8}{\eta}  \mu^t  \sum_{k=1}^n \| x_k(0) - \alpha(0) g_k(0) \|  +  \frac{8nG}{\eta(1-\mu)}   \big(  \alpha(0) \mu^{t/2}  +  \alpha(\ceil{t/2}) \big),
\end{align*}
where $\floor{\cdot}$ denotes the floor function. 
Note that $\| z_i(t+1) - \frac{1}{n} \sum_{j=1}^n (x_j(t) -\alpha(t) g_j(t))\| = \| z_i(t+1) - \frac{ h(t)^\top \1}{n} \|$. The lemma thus has been proved. 
\hfill $\qed$

We are now in a position to first prove Theorem~\ref{thm:bound_everage_n_convex_bound}.

{\bf Proof of Theorem~\ref{thm:bound_everage_n_convex_bound}:}
Note that for all $i \in \mathcal{V}$ and $t \ge 0$, 
\begin{align}
    & \big\|\langle z(t+1) \rangle - z_i(t+1) \big\|+\big\|\bar z(t+1) - z_i(t+1) \big\| \nonumber\\ 
    \le \; & \Big\|\langle z(t+1) \rangle - \frac{1}{n} \sum_{k=1}^n (x_k(t) -\alpha(t) g_k(t)) \Big\| 
    + \Big\|\bar z(t+1) - \frac{1}{n} \sum_{k=1}^n (x_k(t) -\alpha(t) g_k(t)) \Big\| \nonumber\\
    & + 2 \Big\| z_i(t+1) - \frac{1}{n} \sum_{k=1}^n (x_k(t) -\alpha(t) g_k(t)) \Big\| \nonumber\\
    = \; & \Big\|\sum_{j=1}^n \pi_j \big(z_j(t+1)  - \frac{1}{n} \sum_{k=1}^n (x_k(t) -\alpha(t) g_k(t)) \big) \Big\| 
    +  \Big\|\sum_{j=1}^n \frac{1}{n} \big(z_j(t+1)  - \frac{1}{n} \sum_{k=1}^n (x_k(t) -\alpha(t) g_k(t)) \big) \Big\|  \nonumber\\
    & + 2 \Big\| z_i(t+1) - \frac{1}{n} \sum_{k=1}^n (x_k(t) -\alpha(t) g_k(t)) \Big\| \nonumber\\
    \le \; & \sum_{j=1}^n \pi_j(t)  \Big\|  z_j (t+1) - \frac{1}{n} \sum_{k=1}^n (x_k(t) -\alpha(t) g_k(t)) \Big\|  
    + \sum_{j=1}^n \frac{1}{n} \Big\|  z_j (t+1) - \frac{1}{n} \sum_{k=1}^n (x_k(t) -\alpha(t) g_k(t)) \Big\|  \nonumber\\ 
    & + 2 \Big\|  z_i (t+1) - \frac{1}{n} \sum_{k=1}^n (x_k(t) -\alpha(t) g_k(t)) \Big\|  \nonumber\\ 
    \le \; & \frac{32}{\eta}  \mu^t \Big\| \sum_{k=1}^n  x_k(0) - \alpha(0) g_i(0) \Big\| +  \frac{32nG}{\eta}   \sum_{s=0}^t  \mu^{t-s} \alpha(s),  \label{eq:bound_for_difference2_general}
\end{align}
where we used Lemma~\ref{lemma:bound_consensus_push_SA} in the last inequality.
Similarly, if $\{ \alpha(t) \}$ satisfies Assumption~\ref{assum:step-size},
\begin{align}
    &\big\|\langle z(t+1) \rangle - z_i(t+1) \big\|+\big\|\bar z(t+1) - z_i(t+1) \big\| \nonumber\\ 
    \le \;&  \frac{32}{\eta}  \mu^t \Big\| \sum_{k=1}^n  x_k(0) - \alpha(0) g_i(0) \Big\|  +  \frac{32nG}{\eta(1-\mu)}   \Big(  \alpha(0) \mu^{t/2}  +  \alpha(\ceil{t/2}) \Big). \label{eq:bound_for_difference2}
\end{align}
From \eqref{eq:<z>update},
\begin{align}
    \big\| \langle z(t+1) \rangle - z^* \big\|^2
    &= \Big\| \langle z(t) \rangle - z^* - \frac{\alpha(t)}{n}\sum_{i=1}^n g_i(t) \Big\|^2 \nonumber\\
    &\le \big\| \langle z(t) \rangle - z^* \big\|^2 + \Big\| \frac{\alpha(t)}{n}\sum_{i=1}^n g_i(t) \Big\|^2  - 2 \big(\langle z(t) \rangle - z^*\big)^\top \Big(\frac{\alpha(t)}{n}\sum_{i=1}^n g_i(t)\Big) \nonumber \\
    & \le \big\| \langle z(t) \rangle - z^* \big\|^2 + \alpha^2(t) G^2  - 2 \big(\langle z(t) \rangle - z^*\big)^\top \Big(\frac{\alpha(t)}{n}\sum_{i=1}^n g_i(t)\Big). \label{eq:equation 1}
\end{align}
Moreover,
\begin{align}
    \big(\langle z(t) \rangle - z^*\big)^\top g_i(t)
    & = \big(\langle z(t) \rangle - z_i(t) \big)^\top g_i(t) + \big( z_i(t) - z^*\big)^\top g_i(t) \nonumber \\
    & \ge f_i(z_i(t))  - f_i(z^*) - G \big\|\langle z(t) \rangle - z_i(t) \big\|  \label{eq:equation 4_1} \\ 
    & \ge f_i(\bar z(t))  - f_i(z^*) - G \big\|\langle z(t) \rangle - z_i(t) \big\| 
    - G \big\|\bar z(t) - z_i(t) \big\|, \label{eq:equation 4}
\end{align}
where we used \eqref{eq:subgradient} and \eqref{eq:G} in deriving \eqref{eq:equation 4_1}, and used \eqref{eq:G} to get \eqref{eq:equation 4}.
Combining \eqref{eq:equation 1} and \eqref{eq:equation 4}, 
\begin{align*}
    \big\| \langle z(t+1) \rangle - z^* \big\|^2  
    \le \;& \big\| \langle z(t) \rangle - z^* \big\|^2 + \alpha^2(t) G^2 - 2\alpha(t) \big( f(\bar z(t) )  - f(z^*) \big)\\
    & + \frac{2G\alpha(t)}{n} \sum_{i=1}^n \big( \big\|\langle z(t) \rangle - z_i(t) \big\| + \big\|\bar z(t) - z_i(t) \big\| \big),
\end{align*}
which implies that
\begin{align*}
    2\alpha(t) \big( f(\bar z(t) )  - f(z^*) \big)
    \le \;& \big\| \langle z(t) \rangle - z^* \big\|^2 + \alpha^2(t) G^2 - \big\| \langle z(t+1) \rangle - z^* \big\|^2 \\
    & + \frac{2G\alpha(t)}{n} \sum_{i=1}^n \big( \big\|\langle z(t) \rangle - z_i(t) \big\| + \big\|\bar z(t) - z_i(t) \big\| \big).
\end{align*}
Summing this relation over time, it follows that 
\begin{align*}
    \sum_{\tau =0}^t 2\alpha(\tau) \big( f(\bar z(\tau) )  - f(z^*) \big)
    \le \;& \big\| \langle z(0) \rangle - z^* \big\|^2  - \big\| \langle z(t+1) \rangle - z^* \big\|^2 + G^2\sum_{\tau =0}^t \alpha^2(\tau) \\
    & + \sum_{\tau =0}^t  \frac{2G\alpha(\tau)}{n} \sum_{i=1}^n \big( \big\|\langle z(\tau) \rangle - z_i(\tau) \big\| + \big\|\bar z(\tau) - z_i(\tau) \big\| \big).
\end{align*}
Note that
\begin{align*}
    f\bigg(\frac{\sum_{\tau =0}^t \alpha(\tau) \bar z(\tau) }{\sum_{\tau =0}^t \alpha(\tau)}\bigg) - f(z^*)
     \le \frac{ \sum_{\tau =0}^t 2\alpha(\tau) ( f(\bar z(\tau)) - f(z^*) ) }{\sum_{\tau =0}^t 2\alpha(\tau)}.
\end{align*}
It follows that
\begin{align}
    f\bigg(\frac{\sum_{\tau =0}^t \alpha(\tau) \bar z(\tau) }{\sum_{\tau =0}^t \alpha(\tau)}\bigg) - f(z^*)
    \le \; & \frac{ \| \langle z(0) \rangle - z^* \|^2  - \| \langle z(t+1) \rangle - z^* \|^2 + G^2\sum_{\tau =0}^t \alpha^2(\tau) }{\sum_{\tau =0}^t 2\alpha(\tau) }  \nonumber\\
    & + \frac{\sum_{\tau =0}^t \frac{2G\alpha(\tau)}{n} \sum_{i=1}^n (\|\langle z(\tau) \rangle - z_i(\tau) \| + \|\bar z(\tau)  - z_i(\tau) \|)}{\sum_{\tau =0}^t 2\alpha(\tau) } \nonumber\\
    \le \; & \frac{\sum_{\tau =0}^t {G\alpha(\tau)} \sum_{i=1}^n (\|\langle z(\tau) \rangle - z_i(\tau) \| + \|\bar z(\tau)  - z_i(\tau) \|)}{{n}\sum_{\tau =0}^t \alpha(\tau) }  \nonumber\\
    & +  \frac{ \| \langle z(0) \rangle - z^* \|^2 + G^2\sum_{\tau =0}^t \alpha^2(\tau) }{\sum_{\tau =0}^t 2\alpha(\tau) } \label{eq:bound_mid}.
\end{align} 

We next consider the time-varying and fixed stepsizes separately. 

1) If the stepsize $\alpha(t)$ is time-varying and satisfies Assumption~\ref{assum:step-size}, then combining \eqref{eq:bound_for_difference2} and \eqref{eq:bound_mid}, 
\begin{align*}
    f\bigg(\frac{\sum_{\tau =0}^t \alpha(\tau) \bar z(\tau) }{\sum_{\tau =0}^t \alpha(\tau)}\bigg) - f(z^*) 
    \le \;& \frac{ \| \langle z(0) \rangle - z^* \|^2 + G^2\sum_{\tau =0}^t \alpha^2(\tau) }{\sum_{\tau =0}^t 2\alpha(\tau) }  \\
    &  + \frac{ G\alpha(0) \sum_{i=1}^n (\|\langle z(0) \rangle - z_i(0) \| + \|\bar z(0)  - z_i(0) \|)}{n \sum_{\tau =0}^t \alpha(\tau) } \\
     & + \frac{32G}{\eta} \Big\| \sum_{i=1}^n  x_i(0) - \alpha(0) g_i(0) \Big\| \frac{\sum_{\tau =0}^{t-1} \alpha(\tau) \mu^\tau  }{\sum_{\tau =0}^t \alpha(\tau) } \\
     & + \frac{32nG^2}{\eta(1-\mu)} \cdot\frac{\sum_{\tau =0}^{t-1} \alpha(\tau)   (  \alpha(0) \mu^{\tau/2}  +  \alpha(\ceil{\tau/2}) )}{\sum_{\tau =0}^t \alpha(\tau) }.
\end{align*} 
From Proposition~\ref{lemma:push-sum_pi_intfty} and $y_i(0) = 1$, $\pi_i(0) = \frac{1}{n}$ for all $i\in\mathcal{V}$, which implies that $ \langle z(0) \rangle = \frac{1}{n}\sum_i^n z_i(0) = \bar z(0)$. We thus have derived \eqref{eq:bound_timevarying}.

2) If the stepsize is fixed and  $\alpha(t) = 1/\sqrt{T}$ for all $t\ge0$, then from \eqref{eq:bound_mid}, 
\begin{align*}
    f\bigg(\frac{\sum_{\tau =0}^{T-1} \bar z(\tau)  }{ T } \bigg) - f(z^*)
    \le \;& \frac{{G} \sum_{\tau =0}^{T-1} \sum_{i=1}^n (\|\langle z(\tau) \rangle - z_i(\tau) \|+\|\bar z(\tau) - z_i(\tau) \|)}{{n}T } 
    + \frac{ \| \langle z(0) \rangle - z^* \|^2 + G^2 }{ 2\sqrt{T} }\\
    \overset{(a)}{\le} \; & \frac{G \sum_{i=1}^n (\|\langle z(0) \rangle - z_i(0) \|+\|\bar z(0) - z_i(0) \|)}{nT } 
    + \frac{ \| \langle z(0) \rangle - z^* \|^2 + G^2 }{ 2\sqrt{T} } \\
    &+ \frac{32n G^2 }{T \eta} \sum_{\tau =0}^{T-2} \sum_{s=0}^\tau  \mu^{\tau-s} \frac{1}{\sqrt{T}} 
    + \frac{32 G }{T \eta} \Big\| \sum_{i=1}^n  x_i(0) - \frac{1}{\sqrt{T}} g_i(0) \Big\| \sum_{\tau =0}^{T-2} \mu^\tau 
    \\
    \le \;&  \frac{G \sum_{i=1}^n (\|\langle z(0) \rangle - z_i(0) \|+\|\bar z(0) - z_i(0) \|)}{{n}T }
    +\frac{ \| \langle z(0) \rangle - z^* \|^2 + G^2 }{ 2\sqrt{T} } \\
    &+ \frac{32n G^2 }{ \sqrt{T}\eta(1- \mu)} 
    + \frac{32 G }{T \eta(1 - \mu)} \Big\| \sum_{i=1}^n  x_i(0) - \frac{1}{\sqrt{T}} g_i(0) \Big\|,
\end{align*}
where we used \eqref{eq:bound_for_difference2_general} in (a).
Since $ \langle z(0) \rangle = \frac{1}{n}\sum_i^n z_i(0) = \bar z(0)$, we have derived \eqref{eq:bound_fixed}.
\hfill $\qed$

We next prove Theorem~\ref{thm:bound_everage_n_convex}.

{\bf Proof of Theorem~\ref{thm:bound_everage_n_convex}:}
We consider the time-varying and fixed stepsizes separately. First, if the stepsize $\alpha(t)$ is time-varying and satisfies Assumption~\ref{assum:step-size}, then 
\begin{align*}
    & \lim_{t\to\infty}\frac{ \| \langle z(0) \rangle - z^* \|^2 + G^2\sum_{\tau =0}^t \alpha^2(\tau) }{\sum_{\tau =0}^t 2\alpha(\tau) } = 0,\\
    &\lim_{t\to\infty} \frac{ \sum_{i=1}^n (\|\langle z(0) \rangle - z_i(0) \| + \|\bar z(0)  - z_i(0) \|)}{\sum_{\tau =0}^t \alpha(\tau) } = 0.
\end{align*}
Note that 
\begin{align*}
    \sum_{\tau =0}^{t-1} \alpha(\tau) \mu^\tau \le \frac{\alpha(0)   }{1-\mu}, \;\;\;\;\;
    \sum_{\tau =0}^{t-1} \alpha(\tau)   \big(  \alpha(0) \mu^{\tau/2}  +  \alpha(\ceil{\tau/2}) \big) 
    \le \frac{\alpha^2(0) }{1- \mu^{1/2}}  +  \sum_{\tau =0}^{t-1} \alpha^2(\ceil{\tau/2}).
\end{align*}
It follows that
\begin{align*}
    \lim_{t\to\infty} \frac{\sum_{\tau =0}^{t-1} \alpha(\tau) \mu^\tau  }{\sum_{\tau =0}^t \alpha(\tau) } = 0,\;\;\;\;\;
    \lim_{t\to\infty} \frac{\sum_{\tau =0}^{t-1} \alpha(\tau)   (  \alpha(0) \mu^{\tau/2}  +  \alpha(\ceil{\tau/2}) )}{\sum_{\tau =0}^t \alpha(\tau) } = 0.
\end{align*}
Then, the statement in item 1) follows from \eqref{eq:bound_timevarying}. Second, for the case when the stepsize is fixed and  $\alpha(t) = 1/\sqrt{T}$ for $T>0$ steps, the statement in item 2) is an immediate consequence of \eqref{eq:bound_fixed}.
\hfill $\qed$

We finally prove Theorem~\ref{thm:bound_everage_zi}.

{\bf Proof of Theorem~\ref{thm:bound_everage_zi}:}
The proof makes use of the same logic as the proof of Theorem~\ref{thm:bound_everage_n_convex_bound}. Key technical difference lies in three inequalities. 
First, using similar arguments to
those in deriving \eqref{eq:bound_for_difference2_general}, for all $i,k \in \mathcal{V}$ and $t \ge 0$,
\begin{align}
    & \|\langle z(t+1) \rangle - z_i(t+1) \|+\| z_k(t+1) - z_i(t+1) \| \nonumber\\
    \le \; & \frac{32}{\eta}  \mu^t  \sum_{j=1}^n \| x_j(0) - \alpha(0) g_j(0) \| +  \frac{32nG}{\eta}   \sum_{s=0}^t  \mu^{t-s} \alpha(s).\label{eq:bound_for_difference2_general_zi}
\end{align}
Second, similarly, if the stepsize sequence $\{ \alpha(t) \}$ satisfies Assumption~\ref{assum:step-size}, for all $i,k \in \mathcal{V}$ and $t \ge 0$,
\begin{align}
    & \|\langle z(t+1) \rangle - z_i(t+1) \|+\|z_k(t+1) - z_i(t+1) \| \nonumber\\ 
    \le \; & \frac{32}{\eta}  \mu^t \sum_{j=1}^n \|  x_j(0) - \alpha(0) g_j(0) \| 
    +  \frac{32nG}{\eta(1-\mu)}   \big(  \alpha(0) \mu^{t/2}  +  \alpha(\ceil{t/2}) \big). \label{eq:bound_for_difference2_zi}
\end{align}
Third, using the same argument as in getting \eqref{eq:equation 4}, for all $i,k \in \mathcal{V}$ and $t \ge 0$,
\begin{align}
    (\langle z(t) \rangle - z^*)^\top g_i(t)
    \ge f_i( z_k(t))  - f_i(z^*) - G \|\langle z(t) \rangle - z_i(t) \| - G \| z_k(t) - z_i(t) \|.
    \label{eq:equation 4_zi}
\end{align}
Note that these inequalities \eqref{eq:bound_for_difference2_general_zi}, \eqref{eq:bound_for_difference2_zi}, \eqref{eq:equation 4_zi} are respectively variants of inequalities \eqref{eq:bound_for_difference2_general}, \eqref{eq:bound_for_difference2}, \eqref{eq:equation 4} in the proof of Theorem~\ref{thm:bound_everage_n_convex_bound}, with $\bar z$ being replaced by $z_k$. The remaining proof steps are the same as those of the proof of Theorem~\ref{thm:bound_everage_n_convex_bound}, except for replacing inequalities \eqref{eq:bound_for_difference2_general}, \eqref{eq:bound_for_difference2}, \eqref{eq:equation 4} with inequalities \eqref{eq:bound_for_difference2_general_zi}, \eqref{eq:bound_for_difference2_zi}, \eqref{eq:equation 4_zi}, respectively. 
\hfill $\qed$

\section{Stochastic Gradient-Push}

Stochastic gradient descent (SGD) is a popular and powerful first-order iterative optimization method in machine learning. It is well known that stochastic gradient descent achieves $O(1/t)$ convergence rate for strongly convex functions \cite{rakhlin2012making,agarwal2009information,nemirovskij1983problem}.
Distributed SGD has attracted considerable attention for convex optimization. Notable examples include papers \cite{yuan2016convergence,pu2020asymptotic} which consider time-invariant undirected graphs and thus use a fixed doubly stochastic matrix. The work of \cite{srivastava2011distributed} extends SGD to time-varying undirected graphs using the averaging scheme with a time-varying doubly stochastic matrix sequence. All the above algorithms achieve $O(1/t)$ convergence rate, which matches the single-agent SGD performance.\footnote{ 
SGD has also been used to tackle distributed non-convex optimization problems \cite{swenson2020distributed} including distributed deep learning \cite{Lian2017,assran19a}. This paper focuses on distributed convex optimization.}

For time-varying directed graphs, Nedi\'c and Olshevsky proposed a push-sum based SGD algorithm, termed stochastic gradient-push, for the scenario when only stochastic gradient samples are available and showed an $O(\ln t/t)$ convergence rate for strongly convex functions whose finite-time error bound is ``compromised'' by tweaking an usual convex combination of historical states \cite{nedic2016stochastic}. 
Apparently there is a convergence rate gap from distributed (multi-agent) SGD over time-varying directed graphs to state-of-the-art (single-agent) SGD, which is stated as an open problem in \cite{nedic2016stochastic}.

This subsection theoretically proves that the stochastic gradient-push algorithm studied in \cite{nedic2016stochastic} actually converges at a rate of $O(1/t)$, which is the same as the best achievable rate of the single-agent counterpart and thus optimal, and therefore resolves the open problem stated in \cite{nedic2016stochastic}.

We begin with the description of the push-sum based stochastic gradient algorithm proposed in \cite{nedic2016stochastic}. 
Consider the same multi-agent optimization problem setting, except that each $f_i$ is a strongly convex cost function, implying that $f(z)$ is also strongly convex (see Lemma \ref{lemma:f-stronglyconvex}) and thus has a unique optimal solution 
which is denoted as $z^*$.
Each agent $i$ updates its variables as follows:\footnote{Following \cite{nedic2016stochastic} we set the starting time at $t=1$. This avoids division by zero when we set specific $\alpha(t)$ later.}
\begin{align}
    x_i(t+1) &= \sum_{j\in\scr{N}_i(t)} w_{ij}(t)\Big[x_j(t) - \alpha(t)\tilde g_j(t)\Big], \;\;\;\;\; x_i(1)\in\R^d,\label{eq:spushsub_x}\\  
    y_i(t+1) &= \sum_{j\in\scr{N}_i(t)} w_{ij}(t)y_j(t),\;\;\;\;\; y_i(1)=1\label{eq:spushsub_y},
\end{align}
where $\alpha(t)$ is the stepsize, $w_{ij}(t)$, $j\in\scr{N}_i(t)$, are positive weights satisfying Assumption \ref{assum:weighted matrix}, and 
$\tilde g_j(t)$ is a stochastic gradient of $f_j(z)$ at $x_j(t)/y_j(t)$ whose definition is given below.


Given any point $x\in\R^d$, a stochastic subgradient of $f_i$ is a noisy sample of the gradient of $f_i$ at $x$, written as
$
    \tilde g_i(x) = \nabla f_i(x) + N_i(x)
$,
where $\nabla f_i(x)$ denotes the gradient of $f_i$ at point $x$ and $N_i(x)$ is an independent random vector with zero mean, i.e., $\mathbf{E}[N_i(x)]=\0$ and thus $\mathbf{E}[\tilde g_i(x)] = \nabla f_i(x)$. 


It is easy to show the algorithm given by \eqref{eq:spushsub_x}--\eqref{eq:spushsub_y} is mathematically equivalent to that in \cite{nedic2016stochastic}.
To state the convergence result of the stochastic gradient-push algorithm \eqref{eq:spushsub_x}--\eqref{eq:spushsub_y}, we need the following standard concepts and assumptions.


A function $f:\R^d \to \R$ is $\lambda$-strongly convex with $\lambda>0$ if for all $x,y \in \R^d$, 
\begin{align} \label{eq:def_strongcov}
    f(x) - f(y) \ge g(y)^\top (x-y) + \frac{\lambda}{2} \| x-y\|^2,
\end{align}
where $g(y)$ is any subgradient of $f$ at point $y$. 

\begin{assumption} \label{assum:stronglyconvex}
Each function $f_i$ is $\lambda_i$-strongly convex with $\lambda_i>0$.
\end{assumption}

A differentiable function $f:\R^d \to \R$ is  $\gamma$-smooth with $\gamma>0$ if for all $x,y \in \R^d$, 
\begin{align} \label{eq:def_smooth1}
    \| \nabla f(x) - \nabla f(y) \| \le \gamma \| x-y\|.
\end{align}
It is easy to show that the condition \eqref{eq:def_smooth1} is equivalent to
\begin{align} \label{eq:def_smooth2}
    f(y)-f(x)-\nabla f(x)^\top (y-x) \le \frac{ \gamma}{2} \| y-x \|^2.
\end{align}


\begin{assumption} \label{assum:smoothness}
Each function $f_i$ is $\gamma_i$-smooth with $\gamma_i> 0$.
\end{assumption}



\begin{lemma}\label{lemma:bound_grad}
If Assumptions \ref{assum:stronglyconvex} and \ref{assum:smoothness} hold and $\{ \bbb{G}(t) \}$ is uniformly strongly connected,
then $\|\nabla f_i(z_i(t))\|$ is uniformly bounded for all $i$ and $t$ almost surely.
\end{lemma}

{\bf Proof of Lemma \ref{lemma:bound_grad}:}
From Theorem~2 in \cite{nedic2016stochastic},  there exists a positive number $D$ such that $\sup_t \|z_i(t)\| \le D$ for all $i\in\scr V$ almost surely.
Then,
\begin{align*}
    \| \nabla f_i(z_i(t))\|^2 
    &\le 2\| \nabla f_i(z_i(t)) - \nabla f_i(z_i(1)) \|^2 + 2\| \nabla f_i(z_i(1)) \|^2 \\  
    &\le 2\gamma_i^2 \| z_i(t) -z_i(1) \|^2 + 2\| \nabla f_i(z_i(1)) \|^2 
    \le 4 \gamma_i^2 D^2 +2\| \nabla f_i(z_i(1)) \|^2,
\end{align*}
which completes the proof. 
\hfill$\qed$

The proof implicitly assumes a bounded 
$\nabla f_i(z_i(1))$, which can be naturally guaranteed by picking a bounded initial point at each agent.
More can be said with the following assumption. 

\begin{assumption} \label{assum:bounded_G}
    For each $i\in\scr V$ and all $x\in\R^d$,  $\|N_i(x)\|\le c_i$ almost surely for some finite positive number $c_i$.
\end{assumption}


\begin{lemma} \label{lemma:bound_G}
If Assumptions \ref{assum:stronglyconvex}--\ref{assum:bounded_G} hold and $\{ \bbb{G}(t) \}$ is uniformly strongly connected,
then $\mathbf{E}[\|\tilde g_i(t)\|]$, $\|\nabla f_i(z_i(t))\|$, and $\| \nabla f_i( \langle z(t) \rangle )\|$ are all uniformly bounded above by a positive number $G$ for all $i$ and $t$ almost surely.
\end{lemma}

{\bf Proof of Lemma~\ref{lemma:bound_G}:}
From Theorem~2 in \cite{nedic2016stochastic},  there exists a positive number $D$ such that $\sup_t \|z_i(t)\| \le D$ for all $i\in\scr V$ almost surely. Then, 
\begin{align*}
    \| \tilde g_i(t)\|^2 
    &\le 2\| \tilde g_i(t) - \nabla f_i(z_i(1)) \|^2 + 2\| \nabla f_i(z_i(1)) \|^2 \\
    &= 2\| \nabla f_i(z_i(t)) -\nabla f_i(z_i(1)) + N_i(z_i(t)) \|^2 + 2\| \nabla f_i(z_i(1)) \|^2  \\  
    & \le 4\| \nabla f_i(z_i(t)) -\nabla f_i(z_i(1))\|^2 + 4\| N_i(z_i(t)) \|^2 + 2\| \nabla f_i(z_i(1)) \|^2  \\
    &\le 4\gamma_i^2 \| z_i(t) -z_i(1) \|^2 + 4 c_i^2 + 2\| \nabla f_i(z_i(1)) \|^2 
    \le 8 \gamma_i^2 D^2 + 4 c_i^2 +2\| \nabla f_i(z_i(1)) \|^2.
\end{align*}
Let $G^2\dfb \max_{i\in\scr V} (8 \gamma_i^2 D^2 + 4 c_i^2 +2\| \nabla f_i(z_i(1)) \|^2)$. 
Then, $\mathbf{E}[\|\tilde g_i(t)\|^2] \le G^2$ for all $i\in\scr V$ and $t$, 
which immediately implies that 
$\mathbf{E}[\|\tilde g_i(t)\|]\le G$.
From the proof of Lemma \ref{lemma:bound_grad}, it is easy to see that $\| \nabla f_i(z_i(t))\| \le G$.
Similar to the proof of Lemma \ref{lemma:bound_grad},
\begin{align*}
    \| \nabla f_i( \langle z(t) \rangle )\|^2 
    &\le 2\| \nabla f_i(\langle z(t) \rangle) - \nabla f_i(z_i(1)) \|^2 + 2\| \nabla f_i(z_i(1)) \|^2 \\
    &\le 2\gamma_i^2 \| \langle z(t) \rangle -z_i(1) \|^2 + 2\| \nabla f_i(z_i(1)) \|^2 \\
    &\le 2\gamma_i^2 \sum_{k=1}^n \pi_k(t) \|z_k(t) -z_i(1)\|^2 + 2\| \nabla f_i(z_i(1)) \|^2 \le 8\gamma_i^2 D^2 + 2\| \nabla f_i(z_i(1)) \|^2,
\end{align*}
which implies that 
$\| \nabla f_i( \langle z(t) \rangle )\|\le G$.
\hfill $\qed$

Assumptions \ref{assum:stronglyconvex}--\ref{assum:bounded_G} are assumed to hold throughout this section. The convergence result of the stochastic gradient-push algorithm \eqref{eq:spushsub_x}--\eqref{eq:spushsub_y} is as follows.

\begin{theorem}\label{thm:rate_z_i_z*}
If $\{ \bbb{G}(t) \}$ is uniformly strongly connected by sub-sequences of length $L$ and $\alpha(t)=2/(\bar \lambda t)$ where $\bar \lambda = \frac{1}{n}\sum_{i=1}^n \lambda_i$, 
then with probability $1$,
\begin{align*}
   &\mathbf{E}\big[ \| z_i(t) - z^*\|^2 \big] 
   \le \frac{8CG^2}{\bar\lambda^2 t} +  \frac{128nG}{\eta(1-\mu)\bar \lambda (t-1)} + \frac{32K_1}{\eta}  \mu^{t-2}   +  \frac{64nG}{\eta(1-\mu)\bar \lambda} \mu^{(t-1)/2}, \;\;\; i\in\scr V,\\
   &\mathbf{E}\big[f(z_i(t))-f(z^*) \big]
     \le  \frac{4CG^2 \bar \gamma}{\bar\lambda^2 t} +  \frac{64nG \bar \gamma}{\eta(1-\mu)\bar \lambda (t-1)} + \frac{16\bar \gamma K_1}{\eta}  \mu^{t-2}   +  \frac{32nG\bar \gamma}{\eta(1-\mu)\bar \lambda} \mu^{(t-1)/2}, \;\;\; i\in\scr V,\\
    &\mathbf{E}\Big[ f\Big(\frac{1}{n} \sum_{k=1}^n x_k(t) \Big)-f(z^*) \Big]
    \le \frac{2CG^2\bar \gamma}{\bar\lambda^2 (t+1)},
\end{align*}
where $G$ is a positive constant defined in Lemma \ref{lemma:bound_G}, $\eta$ and $\mu$ are positive constants respectively satisfying $\eta \ge \frac{1}{n^{nL}}$ and $\mu \le (1-\frac{1}{n^{nL}})^{\frac{1}{L}}$, 
$K_1 = \mathbf{E}[ \sum_{k=1}^n \|  x_k(1) + \alpha(1) \tilde g_k(1) \| ]$, $K_2 = \frac{-\mu^{-\frac{1}{\ln \mu}}}{\ln \mu}$, and $C = 4 +  \frac{128 K_1 K_2 \bar \lambda }{G\eta \mu^2} +  \frac{512n  (K_2+1)}{\eta (1-\mu) \mu }$. 
\end{theorem} 


The above theorem shows that the stochastic gradient-push algorithm guarantees that each agent's state and function value both converge to the unique optimal point at a rate of $O(1/t)$, which respectively match the performance of the classic single-agent SGD (see Lemma 1 and Theorem~1 in \cite{rakhlin2012making}). Thus, stochastic gradient-push is of the optimal convergence rate. 

Setting the stepsize $\alpha(t)$ to be $2/(\bar \lambda t)$ is for the purpose of cleaner bound expressions. As long as it has the form of $c/t$ with a positive constant $c$, the rate of convergence will still be $O(1/t)$, and thus the implementation of the stochastic gradient-push algorithm \eqref{eq:spushsub_x}--\eqref{eq:spushsub_y} does not need network-wide information $\bar\lambda$. 

\subsection{Analysis}

In this subsection, we provide the proof of Theorems~\ref{thm:rate_z_i_z*}. 

Recall $z_i(t)=x_i(t)/y_i(t)$ and $\langle z(t) \rangle = \frac{1}{n} \sum_{i=1}^n y_i(t) z_i(t)$ which is a time-varying convex combination of all $z_i(t)$. Using the same arguments as in the analysis of subgradient-push, the stochastic gradient-push algorithm \eqref{eq:spushsub_x}--\eqref{eq:spushsub_y} can also be rewritten in terms of a $z_i$ dynamics, 
$
z_i(t+1) 
= \sum_{j=1}^n s_{ij}(t)[z_j(t)- \alpha(t)
\frac{\tilde g_j(t)}{y_j(t)}]
$,
and then
\begin{align}
    \langle z(t+1) \rangle
    &= \frac{1}{n} \sum_{i=1}^n y_i(t+1) z_i(t+1)
    = \frac{1}{n} \sum_{i=1}^n y_i(t+1)  s_{ij}(t)\bigg[z_j(t)- \alpha(t) \frac{\tilde g_j(t)}{y_j(t)}\bigg]\nonumber\\
    &= \frac{1}{n} \sum_{j=1}^n y_j(t) \bigg[z_j(t)- \alpha(t) \frac{\tilde g_j(t)}{y_j(t)}\bigg]
    =\langle z(t) \rangle - \frac{\alpha(t)}{n}\sum_{i=1}^n \tilde g_i(t)\label{eq:update_<z>}.
\end{align}
Similar to subgradient-push, the remaining analysis logic is as follows. Using the inequality $\mathbf{E}[ \| z_i(t) - z^*\|^2 ]\le 2\mathbf{E}[ \|\langle z(t) \rangle - z^*\|^2 ] + 2\mathbf{E}[ \|\langle z(t) \rangle - z_i(t)\|^2 ]$, the proof of Theorem \ref{thm:rate_z_i_z*} is then to bound $\mathbf{E}[ \|\langle z(t) \rangle - z^*\|^2 ]$ and $\mathbf{E}[ \|\langle z(t) \rangle - z_i(t)\|^2 ]$ separately. 
The term $\mathbf{E}[ \|\langle z(t) \rangle - z_i(t)\|^2 ]$ can be bounded using consensus related techniques and relatively easy to deal with. 
The analysis mainly focuses on bounding the term $\mathbf{E}[ \|\langle z(t) \rangle - z^*\|^2 ]$ which is a time-varying quadratic Lyapunov comparison function because 
$\mathbf{E}[\|\langle z(t) \rangle - z^*\|^2]=\mathbf{E}[\|\frac{1}{n}\sum_{i=1}^n y_i(t)(z_i(t)-z^*)\|^2] =\mathbf{E}[\|\frac{1}{n}\sum_{i=1}^n x_i(t)-z^*\|^2]$.
Specifically, the bound for $\mathbf{E}[ \|\langle z(t) \rangle - z^*\|^2 ]$ is given in the following proposition.

\begin{proposition} \label{lemma:rate_<z>_z*}
If $\{ \bbb{G}(t) \}$ is uniformly strongly connected by sub-sequences of length $L$ and $\alpha(t)=2/(\bar \lambda t)$ where $\bar \lambda = \frac{1}{n}\sum_{i=1}^n \lambda_i$, then
$
   \mathbf{E}[ \|\langle z(t) \rangle - z^*\|^2 ] \le 4CG^2/(\bar\lambda^2 t)
$ 
almost surely, where $C$ and $G$ are positive constants respectively defined in Theorem~\ref{thm:rate_z_i_z*} and  Lemma~\ref{lemma:bound_G}.
\end{proposition}

To prove the proposition, we need the following lemmas.

\begin{lemma} \label{lemma:bound_sgpush} 
    If $\{ \bbb{G}(t) \}$ is uniformly strongly connected by sub-sequences of length $L$,
    then with probability $1$, for all $i \in \mathcal{V}$ and $t \ge 1$, 
 \begin{align*}
    \mathbf{E}\Big[\| z_i(t+1) - \frac{1}{n} \sum_{j=1}^n (x_j(t) -\alpha(t)\tilde g_j(t))\| \Big]
    &\le \frac{8K_1}{\eta}  \mu^{t-1} +  \frac{8nG}{\eta}  \sum_{s=1}^t \mu^{t-s} \alpha(s).
\end{align*}
If, in addition, Assumption~\ref{assum:step-size} holds, 
then with probability $1$,
for all $i \in \mathcal{V}$ and $t \ge 1$,
\begin{align*}
    \mathbf{E}\Big[\| z_i(t+1) - \frac{1}{n} \sum_{j=1}^n (x_j(t) -\alpha(t)\tilde g_j(t))\|\Big] 
    &\le \frac{8K_1}{\eta}  \mu^{t-1}  +  \frac{8nG}{\eta(1-\mu)}   \big(  \alpha(1) \mu^{t/2}  +  \alpha(\ceil{t/2}) \big).
\end{align*}
Here $K_1$, $\eta>0$ and $\mu\in(0,1)$ are defined in Theorem~\ref{thm:rate_z_i_z*}, Lemma~\ref{lemma:y_bound} and \eqref{mu}, respectively.
\end{lemma}


{\bf Proof of Lemma~\ref{lemma:bound_sgpush}:}
The lemma can be proved using the same arguments as in the proof of Lemma~\ref{lemma:bound_consensus_push_SA}, except for replacing subgradient $g$ by stochastic gradient $\tilde g$. Note that an inequality still holds after taking expectation on both sides and $\mathbf{E}[\|\tilde g_i(t)\|]\le G$ by Lemma \ref{lemma:bound_G}.
\hfill $\qed$

\begin{lemma} \label{lemma:f-stronglyconvex}
If each function $f_i$ is $\lambda_i$-strongly convex with $\lambda_i>0$, then $f$ is $(\frac{1}{n} \sum_{i=1}^n\lambda_i)$-strongly convex.
\end{lemma}
{\bf Proof of Lemma~\ref{lemma:f-stronglyconvex}:}
From \eqref{eq:def_strongcov}, for all $x, y\in \R^d$,  
\begin{align} \label{eq:sum_stonglyconvex}
    f(x) - f(y) = \frac{1}{n} \sum_{i=1}^n f_i(x) - f_i(y) \ge \frac{\sum_{i=1}^n \nabla f_i(y)}{n} (x-y) +  \frac{\sum_{i=1}^n\lambda_i}{2n} \| x-y\|^2,
\end{align}
which completes the proof.
\hfill $\qed$

\begin{lemma} \label{lemma:initial}
$
    \mathbf{E} [\| \langle z(1) \rangle-z^* \|^2 ] \le G^2/\bar \lambda^2
$
almost surely.
\end{lemma}

{\bf Proof of Lemma~\ref{lemma:initial}:}
From Lemma~\ref{lemma:f-stronglyconvex}, $f$ is $\bar \lambda$-strongly convex with $\bar \lambda = \frac{1}{n}\sum_{i=1}^n \lambda_i$. 
Then, $f$ has a unique optimal point $z^*$ with $\nabla f(z^*)=0$, and thus
\begin{align*}
    \nabla f(\langle z(1) \rangle)^\top (\langle z(1) \rangle-z^*)= (\nabla f(\langle z(1) \rangle) - \nabla f(z^*))^\top (\langle z(1) \rangle-z^*) \ge \bar \lambda \| \langle z(1) \rangle-z^* \|^2.
\end{align*}
This, by the Cauchy-Schwartz inequality, implies that 
$
    \|\nabla f(\langle z(1) \rangle)\| 
    \ge \bar \lambda \| \langle z(1) \rangle-z^* \|
$.
Note that
\begin{align*}
    \mathbf{E}\Big[\Big\| \frac{1}{n}\sum_{i=1}^n \tilde g_i(\langle z(1) \rangle) \Big\|^2 \Big] 
    &\le \frac{1}{n}\sum_{i=1}^n \mathbf{E}\big[\| \tilde g_i(\langle z(1) \rangle) \|^2\big] \le G^2
\end{align*}
almost surely, where we used the convexity of squared 2-norm and Lemma \ref{lemma:bound_G}, and 
\begin{align*}
    \mathbf{E}\Big[\Big\| \frac{1}{n}\sum_{i=1}^n \tilde g_i(\langle z(1) \rangle) \Big\|^2 \Big] 
    =& \; \mathbf{E}\Big[\Big\| \frac{1}{n}\sum_{i=1}^n \tilde g_i(\langle z(1) \rangle) - \nabla f(\langle z(1) \rangle) \Big\|^2\Big] + \mathbf{E}\big[\| \nabla f(\langle z(1) \rangle) \|^2\big] \\
    & + 2 \mathbf{E}\Big[ \Big( \frac{1}{n}\sum_{i=1}^n \tilde g_i(\langle z(1) \rangle) - \nabla f(\langle z(1) \rangle) \Big)^\top \nabla f(\langle z(1) \rangle)\Big] \\
    =& \; \mathbf{E}\Big[\Big\| \frac{1}{n}\sum_{i=1}^n \tilde g_i(\langle z(1) \rangle) - \nabla f(\langle z(1) \rangle) \Big\|^2\Big] + \mathbf{E}\big[\| \nabla f(\langle z(1) \rangle) \|^2\big] \\
    & + 2 \mathbf{E}\Big[ \Big( \frac{1}{n}\sum_{i=1}^n N_i(\langle z(1) \rangle)  \Big)^\top \nabla f(\langle z(1) \rangle)\Big] \\
    \ge & \; \mathbf{E}\big[\| \nabla f(\langle z(1) \rangle) \|^2\big],
\end{align*}
where we used the fact that $\mathbf{E}[N_i(x)]=\0$ for all $i\in\scr V$ and $x\in\R^d$ in the last inequality. From the preceding discussion, 
$\bar \lambda^2 \mathbf{E} [\| \langle z(1) \rangle-z^* \|^2 ]\le \mathbf{E}[\| \nabla f(\langle z(1) \rangle) \|^2] \le \mathbf{E}[\| \frac{1}{n}\sum_{i=1}^n \tilde g_i(\langle z(1) \rangle) \|^2] \le G^2$
almost surely, which completes the proof. 
\hfill $\qed$

We are now in a position to prove Proposition~\ref{lemma:rate_<z>_z*}.

{\bf Proof of Proposition~\ref{lemma:rate_<z>_z*}:}
From \eqref{eq:sum_stonglyconvex}, since $\frac{1}{n}\sum_{i=1}^n \nabla f_i(z^*)= \nabla f(z^*) = 0$, for all $x \in \R^d$,
\begin{align} \label{eq:stonglyconvex_lowerbdd}
    f(x) - f(z^*) 
    \ge \frac{\sum_{i=1}^n \nabla f_i(z^*)}{n} (x-z^*) + \frac{\sum_{i=1}^n\lambda_i}{2n} \| x-z^*\|^2 
    = \frac{\sum_{i=1}^n\lambda_i}{2n} \| x-z^*\|^2.
\end{align}
From \eqref{eq:update_<z>},
\begin{align}
    & \|\langle z(t+1) \rangle - z^*\|^2 
    = \Big\| \langle z(t) \rangle - \frac{\alpha(t)}{n}\sum_{i=1}^n \tilde g_i(z_i(t)) - z^*\Big\|^2 \nonumber\\
    = \; & \| \langle z(t) \rangle - z^*\|^2 + \frac{\alpha(t)^2}{n^2} \Big\|\sum_{i=1}^n \tilde g_i(z_i(t))\Big\|^2 - \frac{2\alpha(t)}{n} \sum_{i=1}^n (\tilde g_i(z_i(t)))^\top (\langle z(t) \rangle - z^*). \label{eq:sGP_comp1}
\end{align}
Let $\mathcal{F}(t)$ be all the information generated by the stochastic gradient-push algorithm \eqref{eq:spushsub_x}--\eqref{eq:spushsub_y} until time $t$, i.e., $\mathcal{F}(t) = \{ x_i(s), y_i(s), \tilde g_i(s) : i\in\mathcal{V}, s\in\{1,\ldots,t\} \}$.
From \eqref{eq:def_strongcov}, with probability $1$,
\begin{align*}
    \mathbf{E}\big[\tilde g_i(z_i(t))^\top (z_i(t) - z^*) \;| \; \mathcal{F}(t) \big]
    &= \nabla f_i(z_i(t))^\top (z_i(t) - z^*) 
    \ge f_i(z_i(t)) - f_i(z^*)+\frac{\lambda_i}{2}\|z_i(t) - z^*\|^2\\
    & = f_i(z_i(t)) - f_i(\langle z(t) \rangle) + f_i(\langle z(t) \rangle) - f_i(z^*)+\frac{\lambda_i}{2}\|z_i(t) - z^*\|^2 \\
    &\ge -G\|\langle z(t) \rangle - z_i(t)\|+ f_i(\langle z(t) \rangle) - f_i(z^*)+\frac{\lambda_i}{2}\|z_i(t) - z^*\|^2, 
\end{align*}
where we used the fact that $f_i(y)-f_i(x) \ge - \|\nabla f_i(x)\| \|x-y\|$ for differentiable convex function $f_i$ and Lemma \ref{lemma:bound_G} to obtain the last inequality.
Then, for any $i\in\mathcal{V}$ and $t\ge 1$,
\begin{align*}
    \mathbf{E}\big[\tilde g_i(z_i(t))^\top (\langle z(t) \rangle - z^*) \;| \; \mathcal{F}(t) \big]
    &= \mathbf{E}\big[\tilde g_i(z_i(t))^\top (\langle z(t) \rangle - z_i(t)) + \tilde g_i(z_i(t))^\top (z_i(t) - z^*) \;| \; \mathcal{F}(t) \big]\\
    &= \nabla f_i(z_i(t))^\top (\langle z(t) \rangle - z_i(t)) + \mathbf{E}\big[\tilde g_i(z_i(t))^\top (z_i(t) - z^*) \;| \; \mathcal{F}(t) \big]\\
    &\ge -2G\|\langle z(t) \rangle - z_i(t)\|+ f_i(\langle z(t) \rangle) - f_i(z^*)+\frac{\lambda_i}{2}\|z_i(t) - z^*\|^2
\end{align*}
almost surely, 
which implies that
\begin{align}
    &\frac{1}{n}\sum_{i=1}^n \mathbf{E}\big[(\tilde g_i(z_i(t)))^\top (\langle z(t) \rangle - z^*) \;| \; \mathcal{F}(t) \big] \nonumber\\
    \ge \; & -\frac{2G}{n}\sum_{i=1}^n\|\langle z(t) \rangle - z_i(t)\|+ f(\langle z(t) \rangle) - f(z^*)+ \frac{1}{n}\sum_{i=1}^n \frac{\lambda_i}{2}\|z_i(t) - z^*\|^2 \nonumber\\
    \ge \; & -\frac{2G}{n}\sum_{i=1}^n\|\langle z(t) \rangle - z_i(t)\|+ \frac{\sum_{i=1}^n\lambda_i}{2n} \| \langle z(t) \rangle-z^*\|^2 + \frac{1}{n}\sum_{i=1}^n \frac{\lambda_i}{2}\|z_i(t) - z^*\|^2, \label{eq:sGP_comp2}
\end{align}
where we used \eqref{eq:stonglyconvex_lowerbdd} in the last inequality. Note that for $t\ge 1$, 
\begin{align}
    &\;\;\;\;\; \mathbf{E}\big[ \|\langle z(t+1) \rangle - z_i(t+1) \| \big] \nonumber \\
    & \le \mathbf{E}\big[ \|\langle z(t+1) \rangle - \frac{1}{n} \sum_{k=1}^n (x_k(t) -\alpha(t) \tilde g_k(t)) \| \big]
    + \mathbf{E}\big[ \| z_i(t+1) - \frac{1}{n} \sum_{k=1}^n (x_k(t) -\alpha(t) \tilde g_k(t)) \| \big] \nonumber\\
    & \le \sum_{i=1}^n \frac{y_i(t+1)}{n} \mathbf{E}\big[ \| z_i(t+1) - \frac{1}{n} \sum_{k=1}^n (x_k(t) -\alpha(t) \tilde g_k(t)) \| \big]
    + \mathbf{E}\big[ \| z_i(t+1) - \frac{1}{n} \sum_{k=1}^n (x_k(t) -\alpha(t) \tilde g_k(t)) \| \big] \nonumber\\
    & 
    \le \frac{16 K_1}{\eta}  \mu^{t-1} +  \frac{16nG}{\eta}   \sum_{s=1}^t  \mu^{t-s} \alpha(s)  \label{eq:sGP_comp3}
\end{align}
almost surely, where we used Lemma~\ref{lemma:bound_sgpush} in the last inequality. 
Then, with $\bar \lambda = \frac{1}{n}\sum_{i=1}^n \lambda_i$ and $K_1 = \mathbf{E}[ \sum_{k=1}^n \|  x_k(1) + \alpha(1) \tilde g_k(1) \| ]$, and from \eqref{eq:sGP_comp1}--\eqref{eq:sGP_comp3} and Lemma~\ref{lemma:bound_G}, 
it follows that for $t\ge 2$,
\begin{align*}
    \mathbf{E}\big[ \|\langle z(t+1) \rangle - z^*\|^2 \big]
    \le \;& \big(1- \alpha(t) \bar\lambda\big) \mathbf{E}\big[ \| \langle z(t) \rangle - z^*\|^2\big] + \frac{\alpha(t)^2}{n^2} \mathbf{E}\big[ \|\sum_{i=1}^n \tilde g_i(z_i(t)) \|^2 \big]\\
    & +  \frac{4\alpha(t)G}{n}\sum_{i=1}^n \mathbf{E}\big[ \|\langle z(t) \rangle - z_i(t)\|\big] 
    -  \frac{2\alpha(t)}{n}\sum_{i=1}^n \frac{\lambda_i}{2}\mathbf{E}\big[ \|z_i(t) - z^*\|^2\big]\\
    \le \;& \Big(1- \frac{2}{t} \Big)\mathbf{E}\big[ \| \langle z(t) \rangle - z^*\|^2 \big]+ \frac{4G^2}{\bar\lambda^2t^2}   +  \frac{128G K_1}{\eta \bar \lambda t}  \mu^{t-2}  +  \frac{128nG^2}{\eta \bar \lambda t}   \sum_{s=1}^{t-1}  \mu^{t-1-s} \alpha(s)
\end{align*}
with probability $1$. 

To proceed, we claim that with $K_2 = \frac{-\mu^{-\frac{1}{\ln \mu}}}{\ln \mu}$, there holds $\mu^t\le \frac{K_2}{t}$ for any $\mu\in(0,1)$ and $t\ge 1$. To prove this, consider the function $h(t) = t\mu^t$ whose derivative $h'(t) = \mu^t (1 + t \ln \mu) $. 
It is straightforward to verify that $h(t)$ takes the maximum value at $t=\frac{-1}{\ln \mu}$. Thus, $h(t) \le h(\frac{-1}{\ln \mu}) = K_2$, which implies that the claim is true.
Using the similar arguments, it can be shown that $\mu^{\frac{t}{2}}\le \frac{2K_2}{t}$. 
Also note that
\begin{align*}
    \sum_{s=1}^{t-1}  \mu^{t-s} \alpha(s)
    &= \sum_{s=1}^{\lfloor t/2 \rfloor}  \mu^{t-s} \alpha(s) 
    + \sum_{s=\lceil t/2 \rceil}^{t-1}  \mu^{t-s} \alpha(s) 
    \le \alpha(1) \sum_{s=t-\lfloor t/2 \rfloor}^{t-1}  \mu^{s}  
    + \alpha(\lceil t/2 \rceil) \sum_{s=0 }^{t-\lceil t/2 \rceil}  \mu^{s} \\
    &\le \frac{1}{1-\mu}\big( \alpha(1)  \mu^{ t/2 }  
    + \alpha(\lceil t/2 \rceil)  \big ).
\end{align*}
Thus, from the preceding inequalities and $\alpha(t)=2/(\bar \lambda t)$,  with probability $1$ and for $t\ge 2$,
\begin{align}
    &\mathbf{E}\big[ \|\langle z(t+1) \rangle - z^*\|^2 \big] \nonumber\\
    \le \;& \Big(1- \frac{2}{t} \Big)\mathbf{E}\big[ \| \langle z(t) \rangle - z^*\|^2 \big]+ \frac{4G^2}{\bar\lambda^2t^2}   +  \frac{128G K_1 K_2}{\eta \bar \lambda \mu^2t^2} 
    +  \frac{128nG^2}{\eta \bar \lambda  (1-\mu)  \mu t}    \Big(  \alpha(1) \mu^{\frac{t}{2}}  +  \frac{2}{\bar \lambda \ceil{\frac{t}{2}}} \Big) \nonumber\\
    \le \;& \Big(1- \frac{2}{t} \Big)\mathbf{E}\big[ \| \langle z(t) \rangle - z^*\|^2 \big]+ \frac{4G^2}{\bar\lambda^2 t^2}  +  \frac{128G K_1 K_2}{\eta \bar \lambda \mu^2 t^2}  +  \frac{256nG^2}{\eta \bar \lambda^2  (1-\mu) \mu t }    \Big(  \mu^{\frac{t}{2}}  +  \frac{2}{t} \Big) \nonumber\\
    \le \;& \Big(1- \frac{2}{t} \Big) \mathbf{E}\big[ \| \langle z(t) \rangle - z^*\|^2 \big]+ \frac{CG^2}{\bar\lambda^2 t^2},\label{eq:yixuaninduction}
\end{align}
in which $C$ is defined in Theorem \ref{thm:rate_z_i_z*}. 


We next prove the proposition statement 
$
   \mathbf{E}[ \|\langle z(t) \rangle - z^*\|^2 ] \le 4CG^2/(\bar\lambda^2 t)
$
almost surely by induction on $t$. 
To this end, we will treat $t=2$ as the base case, and thus first consider $t=1$ separately. 
When $t=1$, from Lemma~\ref{lemma:initial} and the definition of $C$, it is easy to see that  
$
    \mathbf{E}[ \| \langle z(1) \rangle-z^* \|^2] \le 16G^2/\bar\lambda^2 \le 4CG^2/\bar\lambda^2
$ almost surely, thus in this case the statement is true.

For the base case when $t=2$, from \eqref{eq:update_<z>}, 
\begin{align*}
    \mathbf{E}\big[ \| \langle z(2) \rangle-z^* \|^2\big] 
    &= \mathbf{E}\big[\| \langle z(1) \rangle - z^* - \frac{\alpha(1)}{n}\sum_{i=1}^n \tilde g_i(z_i(1)) \|^2 \big] \nonumber\\
    &\le  2\mathbf{E}\big[  \| \langle z(1) \rangle - z^*\|^2\big] + \frac{2\alpha(1)^2}{n^2}\mathbf{E}\big[\sum_{i=1}^n \|\tilde g_i(z_i(1))\|\big]^2\\
    &\overset{(a)}{\le}  2\mathbf{E}\big[  \| \langle z(1) \rangle - z^*\|^2\big] + \frac{2\alpha(1)^2}{n}\sum_{i=1}^n \mathbf{E}\big[\|\tilde g_i(z_i(1))\|^2\big]
    \overset{(b)}{\le}  \frac{2G^2}{\bar \lambda^2} + \frac{8G^2 }{\bar \lambda^2 } \overset{(c)}{\le} \frac{CG^2}{\bar\lambda^2}
\end{align*}
almost surely, where we used the Cauchy--Schwarz inequality in getting (a),  
used Lemma \ref{lemma:initial} and the fact $\mathbf{E}\|\tilde g_i(z_i(1))\|^2\le G^2$ for all $i$ (see the proof of Lemma~\ref{lemma:bound_G}) to obtain inequality (b), and used the fact that $C \ge  \frac{512n}{\eta (1-\mu) \mu } \ge 512$ in inequality (c). 

Now suppose that the statement holds for $t=\tau$ where $\tau\ge 2$ is a positive integer. 
Then, from \eqref{eq:yixuaninduction},
\begin{align*}
    \mathbf{E}\big[ \|\langle z(\tau+1) \rangle - z^*\|^2 \big] & \le \Big(1- \frac{2}{\tau} \Big) \frac{4CG^2}{\bar\lambda^2 \tau} + \frac{CG^2}{\bar\lambda^2 \tau^2}
    = \frac{4CG^2}{\bar\lambda^2 } \cdot\frac{\tau-\frac{7}{4}}{\tau^2}.
\end{align*}
Since for any $\tau\ge 1$, it is easy to show that 
$
    (\tau-\frac{7}{4})/\tau^2 < 1/(\tau+1)
$,
the above inequality immediately implies that the statement holds for $t=\tau+1$. By induction, the statement is established.
\hfill $\qed$

It is easy to prove Theorem~\ref{thm:rate_z_i_z*} with Proposition~\ref{lemma:rate_<z>_z*}.

{\bf Proof of Theorem~\ref{thm:rate_z_i_z*}:}
It is easy to see that $\{\alpha(t) = 2/(\bar \lambda t)\}$ satisfies Assumption~\ref{assum:step-size}.
Then, from Lemma~\ref{lemma:bound_sgpush}, the inequality \eqref{eq:sGP_comp3} can be further bounded as
\begin{align*}
    \mathbf{E}\big[ \|\langle z(t+1) \rangle - z_i(t+1) \| \big]
    \le \frac{16 K_1}{\eta}  \mu^{t-1}  +  \frac{16nG}{\eta(1-\mu)}   \big(  \alpha(1) \mu^{t/2}  +  \alpha(\ceil{t/2}) \big)
\end{align*}
with probability $1$.
From this inequality and Proposition~\ref{lemma:rate_<z>_z*}, 
\begin{align*}
    \mathbf{E}\big[ \| z_i(t) - z^*\|^2 \big]
    &\le 2\mathbf{E}\big[ \langle z(t) \rangle - z^*\|^2 \big] + 2\mathbf{E}\big[ \langle z(t) \rangle - z_i(t)\|^2 \big] \\
    &\le  \frac{8CG^2}{\bar\lambda^2 t} + \frac{32K_1}{\eta}  \mu^{t-2}  +  \frac{32nG}{\eta(1-\mu)}   \big(  \alpha(1) \mu^{(t-1)/2}  +  \alpha(\ceil{(t-1)/2}) \big) \\
    &\le  \frac{8CG^2}{\bar\lambda^2 t} + \frac{32K_1}{\eta}  \mu^{t-2}  +    \frac{64nG}{\eta(1-\mu)\bar \lambda} \mu^{(t-1)/2} + \frac{128nG}{\eta(1-\mu)\bar \lambda (t-1)}
\end{align*}
almost surely, which proves the first inequality in the theorem statement. 

Since $\nabla f(z^*)=0$ and $f$ is $\bar \lambda$-strongly convex (cf. Lemma \ref{lemma:f-stronglyconvex}), \eqref{eq:def_smooth2} implies that  
$\mathbf{E}[f(z_i(t))-f(z^*) ] 
     \le (\bar \gamma/2) \mathbf{E}[\| z_i(t)-z^* \|^2]$.
Then, the second inequality in the theorem statement immediately follows from the first one. 
Similarly, \eqref{eq:def_smooth2} also implies that $\mathbf{E}[f(\langle z(t) \rangle)-f(z^*) ] 
     \le (\bar \gamma/2) \mathbf{E}[\| \langle z(t) \rangle-z^* \|^2]$,
and thus the last inequality in the theorem statement is a direct consequence of 
Proposition~\ref{lemma:rate_<z>_z*}.
\hfill $\qed$

\section{Heterogeneous Multi-Agent Systems}\label{sec:mix}


In both subgradient-push and stochastic gradient-push, all the agents in a multi-agent network perform the same order of operations, namely an optimization step (subgradient or stochastic gradient descent) followed by the push-sum updates. In this section, we aim to relax this order restriction. To illustrate the idea, we focus on the subgradient method, and expect that the idea also works for stochastic gradient descent and other push-sum based first-order optimization methods.

A simple example is the following variant of subgradient-push in which the order of subgradient and push-sum operations is swapped. To be more precise, each agent $i$ updates its variables as
\begin{align}
    x_i(t+1) &= \sum_{j\in\scr{N}_i(t)} w_{ij}(t)x_j(t) - \alpha(t)g_i(t), \;\;\;\;\; x_i(0)\in\R^d,\label{eq:pushfirst_x}\\  y_i(t+1) &= \sum_{j\in\scr{N}_i(t)} w_{ij}(t)y_j(t),\;\;\;\;\; y_i(0)=1\label{eq:pushfirst_y},
\end{align}
where $\alpha(t)$, $w_{ij}(t)$, and $g_i(t)$ are the same as those in subgradient-push.
In the above algorithm \eqref{eq:pushfirst_x}--\eqref{eq:pushfirst_y} each agent $i$ performs the push-sum updates first for both variables and then the subgradient update for $x_i$ variable. We thus call the algorithm push-subgradient.

Push-subgradient can achieve the same performance as subgradient-push, namely it converges to an optimal solution at a rate of $O(1/\sqrt{t})$ for general convex functions. It turns out that both push-subgradient and subgradient-push are special cases of the following heterogeneous algorithm.

Let $\sigma_i(t)$ be a switching signal of agent $i$ which takes values in $\{0,1\}$. At each time $t$, each agent $j$ transmits two pieces of information, $w_{ij}(t)[x_j(t) - \alpha(t)g_j(t)\sigma_j(t)]$ and $w_{ij}(t)y_j(t)$, to its out-neighbour $i$, and then each agent $i$ updates its variables as follows: 
\begin{align}
    x_i(t+1) &= \sum_{j\in\scr{N}_i(t)} w_{ij}(t)\Big[x_j(t) - \alpha(t)g_j(t)\sigma_j(t)\Big] - \alpha(t)g_i(t)\big(1-\sigma_i(t)\big),\;\;\;\;\; x_i(0)\in\R^d, \label{eq:mix_x}\\  y_i(t+1) &= \sum_{j\in\scr{N}_i(t)} w_{ij}(t)y_j(t),\;\;\;\;\; y_i(0)=1\label{eq:mix_y},
\end{align}
where $\alpha(t)$ is the stepsize, $w_{ij}(t)$, $j\in\scr{N}_i(t)$, are positive weights satisfying Assumption \ref{assum:weighted matrix}, and $g_i(t)$ is a subgradient of $f_i(z)$ at $z_i(t)=x_i(t)/y_i(t)$. 

In the case when all $\sigma_i(t)=1$, $i\in\scr V$, the above algorithm simplifies to the subgradient-push algorithm \eqref{eq:pushsub_x}--\eqref{eq:pushsub_y}. In the case when all $\sigma_i(t)=0$, $i\in\scr V$, the above algorithm simplifies to the push-subgradient algorithm \eqref{eq:pushfirst_x}--\eqref{eq:pushfirst_y}. 
Thus, the heterogeneous algorithm allows each agent to arbitrarily switch between subgradient-push and push-subgradient at any time.

The following theorem shows that the heterogeneous distributed subgradient algorithm \eqref{eq:mix_x}--\eqref{eq:mix_y} still achieves the optimal rate of convergence to an optimal point.

\begin{theorem} \label{thm:bound_everage_n_convex_gene}
    Suppose that $\{ \bbb{G}(t) \}$ is uniformly strongly connected and that $\|g_i(t)\|$ is uniformly bounded for all $i$ and~$t$.
\begin{itemize}
    \item[1)] If the stepsize $\alpha(t)$ is time-varying and satisfies Assumption~\ref{assum:step-size}, then with $z^*\in\scr Z$,
\begin{align*}
    \lim_{t\rightarrow\infty}f\bigg(\frac{\sum_{\tau =0}^t \alpha(\tau) \bar z(\tau) }{\sum_{\tau =0}^t \alpha(\tau)}\bigg) = f(z^*), \;\;\;\;\; 
    \lim_{t\rightarrow\infty}f\bigg(\frac{\sum_{\tau =0}^t \alpha(\tau) z_k(\tau) }{\sum_{\tau =0}^t \alpha(\tau)}\bigg) = f(z^*), \;\;\;\;\; k\in\scr V.
\end{align*}
    
    \item[2)] If the stepsize is fixed and  $\alpha(t) = 1/\sqrt{T}$ for $T>0$ steps, then with $z^*\in\scr Z$,
\begin{align*}
    f\bigg(\frac{\sum_{\tau =0}^{T-1} \bar z(\tau)  }{ T } \bigg) - f(z^*) 
     \le O\Big(\frac{1}{ \sqrt{T}}\Big), \;\;\;\;\;
    f\bigg(\frac{\sum_{\tau =0}^{T-1} z_k(\tau)  }{ T } \bigg) - f(z^*) 
     \le O\Big(\frac{1}{ \sqrt{T}}\Big), \;\;\;\;\; k\in\scr V. 
\end{align*}
\end{itemize}
\end{theorem}

Using similar arguments to those in the proof of Theorem \ref{thm:bound_everage_n_convex}, the above theorem is a consequence of the following theorem.

\begin{theorem} \label{thm:bound_everage_n_convex_bound_gene}
    Suppose that $\{ \bbb{G}(t) \}$ is uniformly strongly connected by sub-sequences of length $L$ and that $\|g_i(t)\|$ is uniformly bounded above by a positive number $G$ for all $i$ and $t$.
\begin{itemize}
    \item[1)] 
If the stepsize $\alpha(t)$ is time-varying and satisfies Assumption~\ref{assum:step-size}, then for all $t \ge 0$,
\begin{align}
    f\bigg(\frac{\sum_{\tau =0}^t \alpha(\tau) \bar z(\tau) }{\sum_{\tau =0}^t \alpha(\tau)}\bigg) - f(z^*) 
    \le \; & \frac{ \| \bar z(0) - z^* \|^2 + G^2\sum_{\tau =0}^t \alpha^2(\tau)}{2\sum_{\tau =0}^t \alpha(\tau) } 
    + \frac{2G\alpha(0) \sum_{i=1}^n \|\bar z(0) - z_i(0) \| }{n\sum_{\tau =0}^t \alpha(\tau) } \nonumber\\
    & + \frac{32G\sum_{i=1}^n \| x_i(0) \|}{\eta} \cdot  \frac{\sum_{\tau =0}^{t-1} \alpha(\tau) \mu^\tau  }{\sum_{\tau =0}^t \alpha(\tau) } \nonumber \\
    & + \frac{32nG^2}{\eta\mu(1-\mu)} \cdot \frac{\sum_{\tau =0}^{t-1} \alpha(\tau)   \big(  \alpha(0) \mu^{\frac{\tau}{2}}  +  \alpha(\ceil{\frac{\tau}{2}}) \big)}{\sum_{\tau =0}^t \alpha(\tau) },\label{eq:bound_timevarying_gene}\\
    f\left(\frac{\sum_{\tau =0}^t \alpha(\tau) z_k(\tau) }{\sum_{\tau =0}^t \alpha(\tau)}\right) - f(z^*) 
    \le \; & \frac{ \| \bar z(0) - z^* \|^2 + G^2\sum_{\tau =0}^t \alpha^2(\tau)}{2\sum_{\tau =0}^t \alpha(\tau) } \nonumber \\
    &+ \frac{G\alpha(0) \sum_{i=1}^n \big(\|\bar z(0) - z_i(0) \| + \|z_k(0)  - z_i(0) \|\big)}{n\sum_{\tau =0}^t \alpha(\tau) } \nonumber\\
    & + \frac{32G\sum_{i=1}^n \| x_i(0) \|}{\eta}  \cdot \frac{\sum_{\tau =0}^{t-1} \alpha(\tau) \mu^\tau  }{\sum_{\tau =0}^t \alpha(\tau) } \nonumber\\
    &+ \frac{32nG^2}{\eta\mu(1-\mu)} \cdot \frac{\sum_{\tau =0}^{t-1} \alpha(\tau)   \big(  \alpha(0) \mu^{\frac{\tau}{2}}  +  \alpha(\ceil{\frac{\tau}{2}}) \big)}{\sum_{\tau =0}^t \alpha(\tau) }, \;k\in\scr V. \label{eq:bound_timevarying_zi_gene}
\end{align} 
\item[2)] If the stepsize is fixed and  $\alpha(t) = 1/\sqrt{T}$ for $T>0$ steps, then
\begin{align}
    f\bigg(\frac{\sum_{\tau =0}^{T-1} \bar z(\tau)  }{ T } \bigg) - f(z^*) 
    \le \; & \frac{ \| \bar z(0) - z^* \|^2 + G^2 }{ 2\sqrt{T} } + \frac{2G \sum_{i=1}^n \|\bar z(0) - z_i(0) \|}{nT } \nonumber\\
    & + \frac{32 G \sum_{i=1}^n \| x_i(0)  \| }{\eta(1 - \mu)T}    + \frac{32nG^2}{ \eta \mu (1- \mu)\sqrt{T}}, \label{eq:bound_fixed_gene}\\
    f\left(\frac{\sum_{\tau =0}^{T-1} z_k(\tau)  }{ T } \right) - f(z^*) 
    \le \; &  \frac{ \| \bar z(0) - z^* \|^2 + G^2 }{ 2\sqrt{T} } + \frac{G \sum_{i=1}^n \big(\|\bar z(0)  - z_i(0) \|+\|z_k(0) - z_i(0) \|\big)}{nT } 
    \nonumber\\
    & + \frac{32 G \sum_{i=1}^n \| x_i(0) \|}{\eta(1 - \mu)T}   + \frac{32nG^2}{ \eta \mu (1- \mu)\sqrt{T}}, \;\; k\in\scr V.  \label{eq:bound_fixed_zi_gene}
\end{align}
\end{itemize}
Here $\eta$ and $\mu$ are positive constants which satisfy $\eta \ge \frac{1}{n^{nL}}$ and $\mu \le (1-\frac{1}{n^{nL}})^{1/L}$, respectively.
\end{theorem}

Theorem \ref{thm:bound_everage_n_convex_bound_gene} is a generalization of Theorems \ref{thm:bound_everage_n_convex_bound} and \ref{thm:bound_everage_zi}, so its proof requires a more complicated treatment than those of Theorems \ref{thm:bound_everage_n_convex_bound} and \ref{thm:bound_everage_zi}. It is not surprising that the bounds in Theorems \ref{thm:bound_everage_n_convex_bound} and~\ref{thm:bound_everage_zi} are slightly better than those in Theorem \ref{thm:bound_everage_n_convex_bound_gene} as the former are tailored for a special case.


\subsection{Analysis}

Same as the subgradient-push and stochastic gradient-push algorithms, the analysis of the heterogeneous distributed subgradient algorithm \eqref{eq:mix_x}--\eqref{eq:mix_y} is built upon the dynamics of $\langle z(t) \rangle$ as follows. 
From Proposition~\ref{lemma:push-sum_pi_intfty}, update \eqref{eq:mix_x}, and Assumption \ref{assum:weighted matrix},
\begin{align} 
    &\langle z(t+1) \rangle = \sum_{i=1}^n \pi_i(t+1) z_i(t+1) = \frac{1}{n}\sum_{i=1}^n x_i(t+1)\nonumber\\
    =\; & \frac{1}{n}\sum_{i=1}^n \sum_{j=1}^n w_{ij}(t)\Big[x_j(t) - \alpha(t)g_j(t)\sigma_j(t)\Big] - \frac{\alpha(t)}{n}\sum_{i=1}^n g_i(t)\big(1-\sigma_i(t)\big) \nonumber\\
    =\; & \frac{1}{n} \sum_{j=1}^n \Big[x_j(t) - \alpha(t)g_j(t)\sigma_j(t)\Big] - \frac{\alpha(t)}{n}\sum_{i=1}^n g_i(t)\big(1-\sigma_i(t)\big) 
     = \langle z(t) \rangle -   \frac{\alpha(t)}{n} \sum_{i=1}^n g_i(t),\label{eq:update_<z>_mixed}
\end{align}
in which $w_{ij}(t)$ is set to zero whenever $j\notin\scr N_i(t)$.

From \eqref{eq:<z>update}, the above dynamics of $\langle z(t) \rangle$ is the same as that in the subgradient-push algorithm \eqref{eq:pushsub_x}--\eqref{eq:pushsub_y}. Similarly, it is also easy to show that the push-subgradient algorithm \eqref{eq:pushfirst_x}--\eqref{eq:pushfirst_y} shares the same $\langle z(t) \rangle$ dynamics. This common dynamics is the basis of the following unified analysis for heterogeneous distributed subgradient.

To prove Theorem \ref{thm:bound_everage_n_convex_bound_gene}, we need the following lemma which is a generalization of Lemma~\ref{lemma:bound_consensus_push_SA}, even though its proof follows the similar flow to that in the proof of Lemma~\ref{lemma:bound_consensus_push_SA}.

 \begin{lemma} \label{lemma:bound_consensus_push_SA_gene} 
    If $\{ \bbb{G}(t) \}$ is uniformly strongly connected by sub-sequences of length $L$ and $\|g_i(t)\|$ is uniformly bounded above by a positive number $G$ for all $i$ and $t$, 
    then for all $t \ge 0$ and $i \in \mathcal{V}$,
\begin{align*}
     \Big\| z_i(t+1) - \frac{1}{n}\sum_{k=1}^n x_k(t) \Big\|  \le \frac{8}{\eta}  \mu^t  \sum_{k=1}^n \| x_k(0) \| + \frac{8 nG}{\eta \mu} \sum_{s=0}^t  \mu^{t-s} \alpha(s).
\end{align*}
If, in addition, Assumption~\ref{assum:step-size} holds, then for all $t \ge 0$ and $i \in \mathcal{V}$,
\begin{align*}
     \Big\| z_i(t+1) - \frac{1}{n}\sum_{k=1}^n x_k(t) \Big\| 
    \le \frac{8}{\eta}  \mu^t  \sum_{k=1}^n \| x_k(0) \| 
    +  \frac{8nG}{\eta \mu (1-\mu)}   \big(  \alpha(0) \mu^{t/2}  +  \alpha(\ceil{t/2}) \big).
\end{align*}
Here $\eta>0$ and $\mu\in(0,1)$ are defined in Lemma~\ref{lemma:y_bound} and \eqref{mu}, respectively.
\end{lemma}

{\bf Proof of Lemma~\ref{lemma:bound_consensus_push_SA_gene}:}
Define $\epsilon_i(t) \dfb \sum_{j=1}^n w_{ij}(t) g_j(t) \sigma_j(t) + g_i(t) (1-\sigma_i(t))$ for each $i\in\scr V$ and 
\begin{align*}
 \epsilon(t)\dfb \begin{bmatrix}\epsilon_1^\top(t)\cr\vdots\cr \epsilon_{n}^\top(t)\end{bmatrix}\in\R^{n\times d}.
\end{align*}
Note that 
\begin{align}
    &\sum_{i=1}^n  \| \epsilon_i(t) \|
    \le  \sum_{i=1}^n \Big(\sum_{j=1}^n w_{ij}(t) \|  g_j(t) \| \sigma_j(t) + \| g_i(t) \|( 1-\sigma_i(t))\Big) \nonumber\\
    &\le G \Big(\sum_{i=1}^n \sum_{j=1}^n w_{ij}(t) \sigma_j(t) + \sum_{i=1}^n (1-\sigma_i(t)) \Big) 
    = G \Big( \sum_{j=1}^n \sigma_j(t) + n - \sum_{i=1}^n \sigma_i(t)) \Big) = nG,\label{eq:bound_sum_epsilon}
\end{align}
in which we used the fact that $\sum_{i=1}^n w_{ij}(t)=1$.
From \eqref{eq:mix_x} and the definition of $\Phi_W(t,\tau)$,
\begin{align*}
    x(t+1) 
    = W(t)x(t) -\alpha(t) \epsilon(t) = \Phi_W(t,0) x(0) - \sum_{l=0}^{t-1} \alpha(l) \Phi_W(t,l+1) \epsilon(l) -\alpha(t) \epsilon(t),
\end{align*}
which implies that
\begin{align}
    W(t+1) x(t+1) &= \Phi_W(t+2,0) x(0) - \sum_{l=0}^{t} \alpha(l) \Phi_W(t+2,l+1) \epsilon(l),\label{eq:h_update_tplus1_gene} \\
    \1^\top x(t+1) 
    & = \1^\top x(0) - \sum_{l=0}^{t} \alpha(l) \1^\top \epsilon(l) \label{eq:h_update_1_gene}.
\end{align}
From Lemma~\ref{lemma:pushsum_product} and \eqref{mu}, there exists a sequence  of stochastic vectors $ \{ \phi(t)\}$ such that for all $i,j \in \mathcal{V}$ and $t\ge s \ge 0$, there holds
$
    | [\Phi_W(t+1,s)]_{ij} - \phi_i(t) |\le 4 \mu^{t-s}
$.
Let $ D(s:t) = \Phi_W(t+1,s) - \phi(t) \1^\top$. 
From \eqref{eq:h_update_tplus1_gene} and \eqref{eq:h_update_1_gene}, 
\begin{align*}
    &W(t+1) x(t+1)  - \phi(t+1)\1^\top x(t+1)\\ 
    =\;& \Phi_W(t+2,0) x(0) - \sum_{l=0}^{t} \alpha(l) \Phi_W(t+2,l+1) \epsilon(l)
    - \phi(t+1) \Big(\1^\top x(0) - \sum_{l=0}^{t}  \alpha(l) \1^\top \epsilon(l)\Big) \\
    =\; & \big(\Phi_W(t+2,0) - \phi(t+1) \1^\top\big) x(0) 
    - \sum_{l=0}^{t} \alpha(l)\big( \Phi_W(t+2,l+1) - \phi(t+1) \1^\top\big) \epsilon(l) \\
    =\; & D(0:t+1) x(0) - \sum_{l=0}^{t} \alpha(l) D(l+1:t+1) \epsilon(l),
\end{align*}
which implies that
\begin{align*}
    x(t+1) = W(t) x(t) - \alpha(t) \epsilon(t) = \phi(t)\1^\top x(t) + D(0:t) x(0) - \sum_{l=0}^{t-1} \alpha(l) D(l+1:t) \epsilon(l) - \alpha(t) \epsilon(t).
\end{align*}
From \eqref{eq:mix_y} and the definition of $D(s:t)$,
$
    y(t+1) = \Phi_W(t+1,0) y(0) = D(0:t) \1 + n \phi(t)
$, or equivalently, $y_i(t+1) = [\Phi_W(t+1,0)\1]_i = [D(0:t) \1]_i + n \phi_i(t)$.
Thus, for all $i\in\scr V$,
\begin{align*}
    & z_i(t+1) - \frac{ x(t)^\top \1}{n} 
    = \frac{x_i(t+1)}{y_i(t+1)} - \frac{x(t)^\top \1}{n}\\
    =\; & \frac{\phi_i(t) x(t)^\top \1 + \sum_{k=1}^n [D(0:t)]_{ik} x_k(0) -\sum_{l=0}^{t-1} \alpha(l) \sum_{k=1}^n [D(l+1:t)]_{ik} \epsilon_k(l) - \alpha(t) \epsilon_i(t)}{[D(0:t) \1]_i + n \phi_i(t)} 
    - \frac{ x(t)^\top \1}{n} \\
    =\; & \frac{n \sum_{k=1}^n [D(0:t)]_{ik} x_k(0)   - [D(0:t) \1]_i x(t)^\top \1 -  n\sum_{l=0}^{t-1} \alpha(l) \sum_{k=1}^n [D(l+1:t)]_{ik} \epsilon_k(l)  - n \alpha(t) \epsilon_i(t) }{n[D(0:t) \1]_i + n^2 \phi_i(t)}.
\end{align*}
From Lemma \ref{lemma:y_bound}, $y_i(t+1) \ge \eta$, so is $[D(0:t) \1]_i + n \phi_i(t)$ for all $i\in \scr V$. 
Thus, 
\begin{align*}
    & \Big\| z_i(t+1) - \frac{ x(t)^\top \1}{n} \Big\| \\
    \le\; & \frac{n \| \sum_{k=1}^n [D(0:t)]_{ik} x_k(0) \|  + \| [D(0:t) \1]_i x(t)^\top \1 \| +   n\sum_{l=0}^{t-1} \alpha(l) \|\sum_{k=1}^n [D(l+1:t)]_{ik} \epsilon_k(l)\|  + n \alpha(t) \|\epsilon_i(t) \|}{n[D(0:t) \1]_i + n^2 \phi_i(t)} \\
    \le\; & \frac{n  (\max_k [D(0:t)]_{ik}) \sum_{k=1}^n  \| x_k(0)\|  + \|[D(0:t) \1]_i x(t)^\top \1 \|  }{n[D(0:t) \1]_i + n^2 \phi_i(t)}\\
    &+ \frac{ n\sum_{l=0}^{t-1} \alpha(l) (\max_k [D(l+1:t)]_{ik}) \sum_{k=1}^n \| \epsilon_k(l)\| + n \alpha(t) \| \epsilon_i(t) \| }{n[D(0:t) \1]_i + n^2 \phi_i(t)}\\
    \le\; & \frac{1}{n\eta} \Big( n (\max_k [D(0:t)]_{ik}) \sum_{k=1}^n  \| x_k(0)\|  + n (\max_k [D(0:t)]_{ik}) \| x(t)^\top \1 \| + n \alpha(t) \| \epsilon_i(t)\| \\
    & \;\;\;\;\;\; +  n\sum_{l=0}^{t-1} \alpha(l) (\max_k [D(l+1:t)]_{ik}) \sum_{k=1}^n \| \epsilon_k(l)\|   \Big) \\
    \le\; &  \frac{1}{ \eta} \Big[  4 \mu^t  \sum_{k=1}^n \| x_k(0)\| 
    + \sum_{l=0}^{t-1} \alpha(l) 4 \mu^{t-l-1}    \sum_{k=1}^n \| \epsilon_k(l)\| + \alpha(t) \| \epsilon_i(t) \| + 4 \mu^t \| x(t)^\top \1 \|\Big]. 
\end{align*}
Also, from \eqref{eq:h_update_1_gene}, 
$
    \| \1^\top x(t+1) \|
     \le \| \1^\top x(0)\| + \| \sum_{l=0}^{t} \alpha(l) \1^\top \epsilon(l)\| 
$.
Then, from the above inequality,
\begin{align*}
    & \Big\| z_i(t+1) - \frac{ x(t)^\top \1}{n} \Big\| \\
    \le\; & \frac{4}{ \eta} \Big[   \mu^{t}  \sum_{k=1}^n \| x_k(0)\| + \sum_{l=0}^{t-1} \alpha(l)  \mu^{t-l-1} \sum_{k=1}^n  \|\epsilon_k(l)\| + \alpha(t) \| \epsilon_i(t) \|  +  \mu^{t} \| \1^\top x(0)\| +  \mu^{t} \Big\| \sum_{l=0}^{t-1} \alpha(l) \1^\top \epsilon(l)\Big\| \Big]\\
    \le\; & \frac{8}{ \eta} \Big[ \mu^{t}  \sum_{k=1}^n \| x_k(0)\| 
    + \sum_{l=0}^{t} \alpha(l)  \mu^{t-l-1} \sum_{k=1}^n \| \epsilon_k(l)\|  \Big]. 
\end{align*}
Using \eqref{eq:bound_sum_epsilon}, it follows that for all $i \in \mathcal{V}$ and $t \ge 0$,
\begin{align*}
    \Big\| z_i(t+1) - \frac{ x(t)^\top \1}{n} \Big\|  \le \frac{8}{\eta}  \mu^t  \sum_{k=1}^n \| x_k(0) \| + \frac{8 nG}{\eta \mu} \sum_{s=0}^t  \mu^{t-s} \alpha(s).
\end{align*}
If the stepsize sequence $\{ \alpha(t) \}$ satisfies Assumption~\ref{assum:step-size}, the above inequality further implies that 
\begin{align*}
    \Big\| z_i(t+1) - \frac{ x(t)^\top \1}{n} \Big\| 
    &\le \frac{8}{\eta}  \mu^t  \sum_{k=1}^n \| x_k(0)  \| 
    +  \frac{8nG}{\eta\mu} \bigg(\sum_{s=0}^{\floor{\frac{t}{2}}}  \mu^{t-s} \alpha(s)
    + \sum_{s=\ceil{\frac{t}{2}}}^{t}  \mu^{t-s} \alpha(s) \bigg)   \\
    & \le \frac{8}{\eta}  \mu^t  \sum_{k=1}^n \| x_k(0) \| 
    +  \frac{8nG}{\eta \mu (1-\mu)}   \big(  \alpha(0) \mu^{t/2}  +  \alpha(\ceil{t/2}) \big).
\end{align*}
This completes the proof.
\hfill $\qed$

We are now in a position to prove Theorem~\ref{thm:bound_everage_n_convex_bound_gene}.

{\bf Proof of Theorem~\ref{thm:bound_everage_n_convex_bound_gene}:}
Note that for all $t \ge 0$ and $i,j\in \mathcal{V}$,
\begin{align}
    &\|\langle z(t+1) \rangle - z_i(t+1) \|+\| z_j(t+1) - z_i(t+1) \| \nonumber\\ 
    \le\;& \Big\|\langle z(t+1) \rangle - \frac{1}{n} \sum_{k=1}^n x_k(t)\Big\| +\Big\|z_j(t+1) - \frac{1}{n} \sum_{k=1}^n x_k(t)\Big\|
    + 2\Big\| z_i(t+1) - \frac{1}{n} \sum_{k=1}^n x_k(t)\Big\|\nonumber\\ 
    \le\;& \sum_{j=1}^n \pi_j(t+1) \Big\|  z_j (t+1) - \frac{1}{n} \sum_{k=1}^n x_k(t) \Big\|+\Big\|z_j(t+1) - \frac{1}{n} \sum_{k=1}^n x_k(t)\Big\|
    + 2\Big\| z_i(t+1) - \frac{1}{n} \sum_{k=1}^n x_k(t)\Big\|\nonumber\\ 
    \le\;& \frac{32}{\eta}  \mu^t  \sum_{i=1}^n \| x_i(0) \| +  \frac{32nG}{\eta \mu}   \sum_{s=0}^t  \mu^{t-s} \alpha(s),\label{eq:bound_for_difference2_general_zi_gene}
\end{align}
where we used Lemma~\ref{lemma:bound_consensus_push_SA_gene} in the last inequality. Similarly, for all $t \ge 0$ and $i\in \mathcal{V}$,
\begin{align}
    &\|\langle z(t+1) \rangle - z_i(t+1) \|+\|\bar z(t+1) - z_i(t+1) \| \nonumber\\ 
    \le\;& \Big\|\langle z(t+1) \rangle - \frac{1}{n} \sum_{k=1}^n x_k(t)\Big\| +\Big\|\bar z(t+1) - \frac{1}{n} \sum_{k=1}^n x_k(t)\Big\|
    + 2\Big\| z_i(t+1) - \frac{1}{n} \sum_{k=1}^n x_k(t)\Big\|\nonumber\\ 
    \le\;& \sum_{j=1}^n \Big(\pi_j(t+1)+\frac{1}{n}\Big) \Big\|  z_j (t+1) - \frac{1}{n} \sum_{k=1}^n x_k(t) \Big\| + 2 \Big\| z_i(t+1) - \frac{1}{n} \sum_{k=1}^n x_k(t)  \Big\|\nonumber\\
    \le\;&   \frac{32}{\eta}  \mu^t  \sum_{k=1}^n \| x_k(0) \| +  \frac{32nG}{\eta\mu}   \sum_{s=0}^t  \mu^{t-s} \alpha(s).\label{eq:bound_for_difference2_general_gene}
\end{align}
If, in addition, the stepsize sequence $\{ \alpha(t) \}$ satisfies Assumption~\ref{assum:step-size}, the above two inequalities can be further bounded by Lemma~\ref{lemma:bound_consensus_push_SA_gene} as follows:
\begin{align}
    \|\langle z(t+1) \rangle - z_i(t+1) \|+\|\bar z(t+1) - z_i(t+1) \|  
    &\le \frac{32}{\eta}  \mu^t  \sum_{k=1}^n \| x_k(0) \|  +  \frac{32nG}{\eta\mu(1-\mu)}   \big(  \alpha(0) \mu^{t/2}  +  \alpha(\ceil{t/2}) \big), \label{eq:bound_for_difference2_gene}\\
    \|\langle z(t+1) \rangle - z_i(t+1) \|+\|z_j(t+1) - z_i(t+1) \| 
    &\le \frac{32}{\eta}  \mu^t \sum_{i=k}^n \|  x_k(0) \| +  \frac{32nG}{\eta\mu(1-\mu)}   \big(  \alpha(0) \mu^{t/2}  +  \alpha(\ceil{t/2}) \big). \label{eq:bound_for_difference2_zi_gene} 
\end{align}
From \eqref{eq:update_<z>_mixed}, for any $z^*\in\scr Z$,
\begin{align}
    \| \langle z(t+1) \rangle - z^* \|^2
    &= \Big\| \langle z(t) \rangle - z^* - \frac{\alpha(t)}{n}\sum_{i=1}^n g_i(t) \Big\|^2 \nonumber\\
    & \le \| \langle z(t) \rangle - z^* \|^2 + \Big\| \frac{\alpha(t)}{n}\sum_{i=1}^n g_i(t) \Big\|^2 - 2 (\langle z(t) \rangle - z^*)^\top \Big(\frac{\alpha(t)}{n}\sum_{i=1}^n g_i(t)\Big) \nonumber \\
    & \le \| \langle z(t) \rangle - z^* \|^2 + \alpha^2(t) G^2  - 2 (\langle z(t) \rangle - z^*)^\top \Big(\frac{\alpha(t)}{n}\sum_{i=1}^n g_i(t)\Big), \label{eq:equation 1_gene}
\end{align}
where we used the convexity of squared 2-norm in the last inequality. 
Moreover, for all $i,k\in\scr V$,
\begin{align}
    (\langle z(t) \rangle - z^*)^\top g_i(t) 
    & = (\langle z(t) \rangle - z_i(t) )^\top g_i(t) + ( z_i(t) - z^*)^\top g_i(t) \nonumber \\
    & \ge f_i(z_i(t))  - f_i(z^*) - G \|\langle z(t) \rangle - z_i(t) \|    \label{eq:equation 4_1_gene} \\
    & \ge f_i( z_k(t))  - f_i(z^*) - G \|\langle z(t) \rangle - z_i(t) \|  - G \| z_k(t) - z_i(t) \|,
    \label{eq:equation 4_zi_gene}
\end{align}
where we used \eqref{eq:subgradient} and \eqref{eq:G} in deriving \eqref{eq:equation 4_1_gene}, and made use of \eqref{eq:G} to get \eqref{eq:equation 4_zi_gene}. 
Similarly, for all $i\in\scr V$,
\begin{align}
     (\langle z(t) \rangle - z^*)^\top g_i(t)
    & \ge f_i(\bar z(t))  - f_i(z^*) - G \|\langle z(t) \rangle - z_i(t) \|  - G \|\bar z(t) - z_i(t) \|. \label{eq:equation 4_gene}
\end{align}
Combining \eqref{eq:equation 1_gene} and \eqref{eq:equation 4_gene}, 
\begin{align*}
    \| \langle z(t+1) \rangle - z^* \|^2  
    \le\; & \| \langle z(t) \rangle - z^* \|^2 + \alpha^2(t) G^2 - 2\alpha(t) ( f(\bar z(t) )  - f(z^*) )\\
    & + \frac{2G\alpha(t)}{n} \sum_{i=1}^n \big( \|\langle z(t) \rangle - z_i(t) \| + \|\bar z(t) - z_i(t) \| \big),
\end{align*}
which implies that
\begin{align*}
    2\alpha(t) ( f(\bar z(t) )  - f(z^*))
    \le \;& \| \langle z(t) \rangle - z^* \|^2 + \alpha^2(t) G^2 - \| \langle z(t+1) \rangle - z^* \|^2 \\
    & + \frac{2G\alpha(t)}{n} \sum_{i=1}^n \big(\|\langle z(t) \rangle - z_i(t) \| + \|\bar z(t) - z_i(t) \| \big).
\end{align*}
Summing this relation over time, it follows that 
\begin{align*}
    \sum_{\tau =0}^t 2\alpha(\tau) ( f(\bar z(\tau) )  - f(z^*) )
    \le \;& \| \langle z(0) \rangle - z^* \|^2  - \| \langle z(t+1) \rangle - z^* \|^2 + G^2\sum_{\tau =0}^t \alpha^2(\tau) \\
    & + \sum_{\tau =0}^t  \frac{2G\alpha(\tau)}{n} \sum_{i=1}^n \big( \|\langle z(\tau) \rangle - z_i(\tau) \| + \|\bar z(\tau) - z_i(\tau) \| \big).
\end{align*}
Then, 
\begin{align}
    f\bigg(\frac{\sum_{\tau =0}^t \alpha(\tau) \bar z(\tau) }{\sum_{\tau =0}^t \alpha(\tau)}\bigg) - f(z^*)
    \le\;& \frac{ \sum_{\tau =0}^t 2\alpha(\tau) ( f(\bar z(\tau)) - f(z^*) ) }{\sum_{\tau =0}^t 2\alpha(\tau)}\nonumber\\
    \le\;& \frac{ \| \langle z(0) \rangle - z^* \|^2  - \| \langle z(t+1) \rangle - z^* \|^2 + G^2\sum_{\tau =0}^t \alpha^2(\tau) }{\sum_{\tau =0}^t 2\alpha(\tau) }  \nonumber\\
    & + \frac{\sum_{\tau =0}^t \frac{2G\alpha(\tau)}{n} \sum_{i=1}^n (\|\langle z(\tau) \rangle - z_i(\tau) \| + \|\bar z(\tau)  - z_i(\tau) \|)}{\sum_{\tau =0}^t 2\alpha(\tau) } \nonumber\\
    \le\;&  \frac{\sum_{\tau =0}^t {G\alpha(\tau)} \sum_{i=1}^n (\|\langle z(\tau) \rangle - z_i(\tau) \| + \|\bar z(\tau)  - z_i(\tau) \|)}{{n}\sum_{\tau =0}^t \alpha(\tau) }  \nonumber\\
    & +  \frac{ \| \langle z(0) \rangle - z^* \|^2 + G^2\sum_{\tau =0}^t \alpha^2(\tau) }{\sum_{\tau =0}^t 2\alpha(\tau) }. \label{eq:bound_mid_gene}
\end{align} 
Similarly, combining \eqref{eq:equation 1_gene} and \eqref{eq:equation 4_zi_gene}, for any $k\in\scr V$, 
\begin{align*}
    \| \langle z(t+1) \rangle - z^* \|^2 
    \le\;& \| \langle z(t) \rangle - z^* \|^2 + \alpha^2(t) G^2 - 2\alpha(t) ( f(z_k(t) )  - f(z^*) )\\
    &  + \frac{2G\alpha(t)}{n} \sum_{i=1}^n \big( \|\langle z(t) \rangle - z_i(t) \| + \| z_k(t) - z_i(t) \| \big),
\end{align*}
which, using the preceding argument, leads to
\begin{align}
    f\bigg(\frac{\sum_{\tau =0}^t \alpha(\tau) z_k(\tau) }{\sum_{\tau =0}^t \alpha(\tau)}\bigg) - f(z^*)
    \le\;&  \frac{\sum_{\tau =0}^t G\alpha(\tau) \sum_{i=1}^n (\|\langle z(\tau) \rangle - z_i(\tau) \| + \|z_k(\tau)  - z_i(\tau) \|)}{{n}\sum_{\tau =0}^t \alpha(\tau) }  \nonumber\\
    & +  \frac{ \| \langle z(0) \rangle - z^* \|^2 + G^2\sum_{\tau =0}^t \alpha^2(\tau) }{\sum_{\tau =0}^t 2\alpha(\tau) }. \label{eq:bound_mid_zi_gene}
\end{align} 

We next consider the time-varying and fixed stepsizes separately. 

1) If the stepsize $\alpha(t)$ is time-varying and satisfies Assumption~\ref{assum:step-size}, then combining \eqref{eq:bound_for_difference2_gene} and \eqref{eq:bound_mid_gene}, 
\begin{align*}
    & f\bigg(\frac{\sum_{\tau =0}^t \alpha(\tau) \bar z(\tau) }{\sum_{\tau =0}^t \alpha(\tau)}\bigg) - f(z^*) \\
    \le\;& \frac{ \| \langle z(0) \rangle - z^* \|^2 + G^2\sum_{\tau =0}^t \alpha^2(\tau) }{\sum_{\tau =0}^t 2\alpha(\tau) }  
    + \frac{ G\alpha(0) \sum_{i=1}^n (\|\langle z(0) \rangle - z_i(0) \| + \|\bar z(0)  - z_i(0) \|)}{n \sum_{\tau =0}^t \alpha(\tau) } \\
     & + \frac{32G}{\eta}  \Big(\sum_{i=1}^n \| x_i(0) \|\Big) \frac{\sum_{\tau =0}^{t-1} \alpha(\tau) \mu^\tau  }{\sum_{\tau =0}^t \alpha(\tau) } 
    + \frac{32nG^2}{\eta\mu(1-\mu)} \cdot \frac{\sum_{\tau =0}^{t-1} \alpha(\tau)   (  \alpha(0) \mu^{\frac{\tau}{2}}  +  \alpha(\ceil{\frac{\tau}{2}}) )}{\sum_{\tau =0}^t \alpha(\tau) }.
\end{align*} 
Similarly, combining \eqref{eq:bound_for_difference2_zi_gene} and \eqref{eq:bound_mid_zi_gene}, 
\begin{align*}
    & f\bigg(\frac{\sum_{\tau =0}^t \alpha(\tau) z_k(\tau) }{\sum_{\tau =0}^t \alpha(\tau)}\bigg) - f(z^*) \\
    \le\;& \frac{ \| \langle z(0) \rangle - z^* \|^2 + G^2\sum_{\tau =0}^t \alpha^2(\tau) }{\sum_{\tau =0}^t 2\alpha(\tau) }  
    + \frac{ G\alpha(0) \sum_{i=1}^n (\|\langle z(0) \rangle - z_i(0) \| + \|z_k(0)  - z_i(0) \|)}{n\sum_{\tau =0}^t \alpha(\tau) } \\
    & + \frac{32G}{\eta}  \Big(\sum_{i=1}^n \| x_i(0) \|\Big) \frac{\sum_{\tau =0}^{t-1} \alpha(\tau) \mu^\tau  }{\sum_{\tau =0}^t \alpha(\tau) } 
    + \frac{32nG^2}{\eta\mu(1-\mu)} \cdot \frac{\sum_{\tau =0}^{t-1} \alpha(\tau)   (  \alpha(0) \mu^{\frac{\tau}{2}}  +  \alpha(\ceil{\frac{\tau}{2}}) )}{\sum_{\tau =0}^t \alpha(\tau) }.
\end{align*} 
For all $i\in\mathcal{V}$, from Proposition~\ref{lemma:push-sum_pi_intfty} and $y_i(0) = 1$, we have $\pi_i(0) = \frac{1}{n}$, which implies that $ \langle z(0) \rangle = \frac{1}{n}\sum_i^n z_i(0) = \bar z(0)$. We thus have derived \eqref{eq:bound_timevarying_gene} and \eqref{eq:bound_timevarying_zi_gene}.

2) If the stepsize is fixed and  $\alpha(t) = 1/\sqrt{T}$ for all $t\ge0$, then from \eqref{eq:bound_mid_gene} and \eqref{eq:bound_for_difference2_general_gene},  
\begin{align*}
    f\bigg(\frac{\sum_{\tau =0}^{T-1} \bar z(\tau)  }{ T } \bigg) - f(z^*)
    \le\;& \frac{{G} \sum_{\tau =0}^{T-1} \sum_{i=1}^n \|\langle z(\tau) \rangle - z_i(\tau) \|+\|\bar z(\tau) - z_i(\tau) \|}{{n}T }
    + \frac{ \| \langle z(0) \rangle - z^* \|^2 + G^2 }{ 2\sqrt{T} } \\
    \le\;& \frac{G \sum_{i=1}^n \|\langle z(0) \rangle - z_i(0) \|+\|\bar z(0) - z_i(0) \|}{nT } + \frac{ \| \langle z(0) \rangle - z^* \|^2 + G^2 }{ 2\sqrt{T} }\\ 
    &  + \frac{32 G }{T \eta}  \Big(\sum_{i=1}^n \| x_i(0)\|\Big) \sum_{\tau =0}^{T-2} \mu^\tau  + \frac{32n G^2 }{T \eta \mu} \sum_{\tau =0}^{T-2} \sum_{s=0}^\tau  \mu^{\tau-s} \frac{1}{\sqrt{T}}\\
    \le\;&  \frac{G \sum_{i=1}^n \|\langle z(0) \rangle - z_i(0) \|+\|\bar z(0) - z_i(0) \|}{{n}T } +\frac{ \| \langle z(0) \rangle - z^* \|^2 + G^2 }{ 2\sqrt{T} } \\
    & + \frac{32 G }{T \eta(1 - \mu)} \sum_{i=1}^n \|  x_i(0) \| + \frac{32n G^2 }{ \sqrt{T}\eta \mu (1- \mu)} .
\end{align*}
Similarly, from \eqref{eq:bound_mid_zi_gene} and \eqref{eq:bound_for_difference2_general_zi_gene}, 
\begin{align*}
    f\bigg(\frac{\sum_{\tau =0}^{T-1} z_k(\tau)  }{ T } \bigg) - f(z^*)
    \le\;& \frac{ G \sum_{\tau =0}^{T-1} \sum_{i=1}^n \|\langle z(\tau) \rangle - z_i(\tau) \|+\|z_k(\tau) - z_i(\tau) \|}{nT }
    + \frac{ \| \langle z(0) \rangle - z^* \|^2 + G^2 }{ 2\sqrt{T} } \\
    \le\;& \frac{G \sum_{i=1}^n \|\langle z(0) \rangle - z_i(0) \|+\|z_k(0) - z_i(0) \|}{nT } + \frac{ \| \langle z(0) \rangle - z^* \|^2 + G^2 }{ 2\sqrt{T} } \\ 
    & + \frac{32 G }{T \eta}  \Big(\sum_{i=1}^n \| x_i(0) \|\Big) \sum_{\tau =0}^{T-2} \mu^\tau + \frac{32n G^2 }{T \eta \mu} \sum_{\tau =0}^{T-2}  \sum_{s=0}^\tau  \mu^{\tau-s} \frac{1}{\sqrt{T}} \\
    \le\;&  \frac{G \sum_{i=1}^n \|\langle z(0) \rangle - z_i(0) \|+\|z_k(0) - z_i(0) \|}{nT }
    +\frac{ \| \langle z(0) \rangle - z^* \|^2 + G^2 }{ 2\sqrt{T} } \\
    & + \frac{32 G }{T \eta(1 - \mu)}  \sum_{i=1}^n \| x_i(0) \|  
    + \frac{32n G^2 }{ \sqrt{T}\eta \mu (1- \mu)}.
\end{align*}
Since $ \langle z(0) \rangle = \frac{1}{n}\sum_i^n z_i(0) = \bar z(0)$, we have derived \eqref{eq:bound_fixed_gene} and \eqref{eq:bound_fixed_zi_gene}.
\hfill $\qed$

\begin{remark}
Theorem \ref{thm:bound_everage_n_convex_bound_gene} and its proof involve a positive constant $\mu$ stemming from \eqref{mu}. Technically, it is possible that $\mu$ in \eqref{mu} equals zero, for example, each $W(t)$ is a rank one column stochastic matrix of the form $u(t)\1^\top$ with $u(t)$ being a column stochastic vector. In this special case, it is easy to see that Theorem \ref{thm:bound_everage_n_convex_bound_gene} and its proof hold for any positive $\mu<1$ because so does inequality \eqref{mu}, yet the bounds in Theorem \ref{thm:bound_everage_n_convex_bound_gene} can be far from being tight. With these in mind, we derive tighter bounds for the special case when $\mu=0$ in this remark. To this end, we first follow the proof of Lemma \ref{lemma:bound_consensus_push_SA_gene}. Note that with $\mu$ being zero $D(s:t)=0$ for all $t\ge s\ge 0$. It follows that $x(t+1) = \phi(t) \1^\top x(t) - \alpha(t) \epsilon(t)$ and $y(t+1) = n\phi(t)$ for all $t\ge 0$ where $ \{ \phi(t)\}$ is a sequence  of stochastic vectors. Then, 
\begin{align*}
    z_i(t+1) - \frac{ x(t)^\top \1}{n} 
    = \frac{\phi_i(t) x(t)^\top \1 - \alpha(t) \epsilon_i(t)}{n \phi_i(t)} 
    - \frac{ x(t)^\top \1}{n} 
    = - \frac{\alpha(t) \epsilon_i(t)}{n \phi_i(t)} 
    = - \frac{\alpha(t) \epsilon_i(t)}{y_i(t+1)}
\end{align*}
for all $i\in\scr V$, which implies that
$
    \| z_i(t+1) - \frac{ x(t)^\top \1}{n} \| 
    \le  \frac{\alpha(t) \| \epsilon_i(t) \|}{y_i(t+1)} \le \frac{\alpha(t)nG}{\eta}
$
where we used Lemma~\ref{lemma:y_bound} and the fact that $\|\epsilon_i(t)\|\le nG$ due to \eqref{eq:bound_sum_epsilon}.
We next follow the proof of Theorem \ref{thm:bound_everage_n_convex_bound_gene}. Using the inequality just derived and the same argument as in the proof of Theorem \ref{thm:bound_everage_n_convex_bound_gene}, 
$\|\langle z(t+1) \rangle - z_i(t+1) \|+\| z_j(t+1) - z_i(t+1) \| \le \frac{4nG\alpha(t)}{\eta}$ and 
$\|\langle z(t+1) \rangle - z_i(t+1) \|+\|\bar z(t+1) - z_i(t+1) \| 
    \le   \frac{4nG\alpha(t)}{\eta}$ for all $t \ge 0$ and $i,j\in \mathcal{V}$. 
Both inequalities are then used in time-varying and fixed stepsize cases. 


If the stepsize $\alpha(t)$ is time-varying and satisfies Assumption~\ref{assum:step-size}, then respectively from \eqref{eq:bound_mid_gene} and \eqref{eq:bound_mid_zi_gene}, 
\begin{align*}
    f\bigg(\frac{\sum_{\tau =0}^t \alpha(\tau) \bar z(\tau) }{\sum_{\tau =0}^t \alpha(\tau)}\bigg) - f(z^*) 
    \le\;& \frac{ \| \bar z(0)  - z^* \|^2 + G^2\sum_{\tau =0}^t \alpha^2(\tau) }{2\sum_{\tau =0}^t \alpha(\tau) } + 
    \frac{ {2G\alpha(0)} \sum_{i=1}^n \|\bar z(0)  - z_i(0) \|}{{n}\sum_{\tau =0}^t \alpha(\tau) } \nonumber\\
    & + \frac{4nG^2}{\eta} \frac{\sum_{\tau =1}^t {\alpha(\tau)}  \alpha(\tau-1)}{\sum_{\tau =0}^t \alpha(\tau) }, \nonumber\\
    f\bigg(\frac{\sum_{\tau =0}^t \alpha(\tau) z_k(\tau) }{\sum_{\tau =0}^t \alpha(\tau)}\bigg) - f(z^*) 
    \le\;& \frac{ \| \bar z(0)  - z^* \|^2 + G^2\sum_{\tau =0}^t \alpha^2(\tau) }{2\sum_{\tau =0}^t \alpha(\tau) } +
    \frac{ G\alpha(0) \sum_{i=1}^n \|\bar z(0)  - z_i(0) \| }{{n}\sum_{\tau =0}^t \alpha(\tau) }  \nonumber\\
    & +   \frac{4nG^2}{\eta} \frac{\sum_{\tau =1}^t \alpha(\tau)  \alpha(\tau - 1)}{\sum_{\tau =0}^t \alpha(\tau) } +
        \frac{ G\alpha(0) \sum_{i=1}^n  \|z_k(0)  - z_i(0) \|}{{n}\sum_{\tau =0}^t \alpha(\tau) }. 
\end{align*} 
If the stepsize is fixed and  $\alpha(t) = 1/\sqrt{T}$ for all $t\ge0$, then respectively from \eqref{eq:bound_mid_gene} and \eqref{eq:bound_mid_zi_gene}, 
\begin{align*}
    f\bigg(\frac{\sum_{\tau =0}^{T-1} \bar z(\tau)  }{ T } \bigg) - f(z^*)
    \le\;&  \frac{2G \sum_{i=1}^n \|\bar z(0) - z_i(0) \|}{{n}T } +\frac{ \| \bar z(0) - z^* \|^2 + G^2 }{ 2\sqrt{T} }  + \frac{4nG^2}{\eta \sqrt{T}}, \\
    f\bigg(\frac{\sum_{\tau =0}^{T-1} z_k(\tau)  }{ T } \bigg) - f(z^*)
    \le\;&  \frac{G \sum_{i=1}^n (\|\bar z(0) - z_i(0) \|+\|z_k(0) - z_i(0) \|)}{nT }
    +\frac{ \| \bar z(0) - z^* \|^2 + G^2 }{ 2\sqrt{T} }  + \frac{4nG^2}{\eta \sqrt{T}}.
\end{align*}
It is clear that all the above bounds are simpler and tighter than those in Theorem \ref{thm:bound_everage_n_convex_bound_gene}.
\hfill$\Box$
\end{remark}

\subsection{A Special Case}

In this subsection, we discuss a special case in which $W(t)$ is a doubly stochastic matrix at all time $t\ge 0$. In this case, it is easy to see from \eqref{eq:mix_y} that $y_i(t)=1$ for all $i\in\scr V$ and $t\ge 0$, and thus $z_i(t)=x_i(t)$ for all $i\in\scr V$ and $t\ge 0$. This observation holds for all push-sum based distributed optimization algorithms studied in this paper as they share the same $y_i(t)$ dynamics which is independent of their $x_i(t)$ dynamics. Then, the subgradient-push, push-subgradient, and heterogeneous subgradient algorithms all simplify to average consensus based subgradient algorithms. Specifically, subgradient-push \eqref{eq:pushsub_x}--\eqref{eq:pushsub_y} simplifies to 
\begin{align}
    x_i(t+1) = \sum_{j\in\scr{N}_i(t)} w_{ij}(t)\Big[x_j(t) - \alpha(t)g_j(x_j(t))\Big], \label{eq:adapt_diffusion}
\end{align}
and push-subgradient \eqref{eq:pushfirst_x}--\eqref{eq:pushfirst_y} simplifies to 
\begin{align}
x_i(t+1) = \sum_{j\in\scr{N}_i(t)} w_{ij}(t)x_j(t) - \alpha(t)g_i(x_i(t)), \label{eq:diffusion_adapt}
\end{align}
which is the ``standard'' average consensus based distributed subgradient proposed in \cite{nedic2009distributed}.
The two updates \eqref{eq:adapt_diffusion} and \eqref{eq:diffusion_adapt} are analogous to the so-called ``adapt-then-combine'' and ``combine-then-adapt'' diffusion strategies in distributed optimization and learning \cite{sayed_magazine}. Thus, in the special case under consideration, the heterogeneous subgradient \eqref{eq:mix_x}--\eqref{eq:mix_y} simplifies to 
\begin{align*}
    x_i(t+1) = \sum_{j\in\scr{N}_i(t)} w_{ij}(t)\Big[x_j(t) - \alpha(t)g_j(x_j(t))\sigma_j(t)\Big] - \alpha(t)g_i(x_i(t))\big(1-\sigma_i(t)\big),
\end{align*}
which is an average consensus based heterogeneous distributed subgradient algorithm allowing each agent to arbitrarily switch between updates \eqref{eq:adapt_diffusion} and \eqref{eq:diffusion_adapt}.
The preceding discussion implies that all the results in this paper apply to the corresponding average consensus based algorithms.

\section{Arbitrary Convex Combination}

In this section, we generalize the push-sum algorithm to solve arbitrary convex combination problems, subsuming the straight average as a special case. 

To this end, let each agent $i$ be aware of a positive number $c_i$ such that $\sum_{i=1}^n c_i = 1$. Each agent $i$ iteratively updates its variables $x_i$ and $y_i$ as follows:
\begin{align} 
    x_i(t+1) &= \sum_{j\in\scr{N}_i(t)} w_{ij}(t)x_j(t),\;\;\;\;\; x_i(0)=c_i x_i^{{\rm int}}\in\R^d,\label{eq:convexx}\\ 
    y_i(t+1) &= \sum_{j\in\scr{N}_i(t)} w_{ij}(t)y_j(t),\;\;\;\;\; y_i(0)=c_i,\label{eq:convexy}
\end{align}
%
where weights $w_{ij}(t)$ are the same as those in \eqref{pushsumx}--\eqref{pushsumy}, so satisfying Assumption \ref{assum:weighted matrix}, 
and $x_i^{{\rm int}}$ denotes the initial value of agent $i$.
Thus, the above algorithm is the same as the standard push-sum except for the initial values of $x_i$ and $y_i$ variables. 
The limiting behavior of the above algorithm is described by the following theorem.

\begin{theorem}\label{thm:convex}
    If $\{ \bbb{G}(t) \}$ is uniformly strongly connected, then each $x_i(t)/y_i(t)$ converges to $\sum_{k=1}^n c_k x_k^{{\rm int}}$ exponentially fast as $t\rightarrow\infty$. 
\end{theorem}

The theorem can be proved using similar arguments to those in the proof of Theorem~\ref{thm:pushsum}, both the conventional proof in the literature and our new proof given in Section \ref{sec:pushsum}. It shows that all the agents can finally figure out an arbitrary given convex combination of their initial values in a distributed manner, provided each agent knows its own convex combination coefficient. 
It is worth noting that although the above algorithm can solve the straight average case, it requires each agent know the value of $1/n$, which is more restrictive than the standard push-sum algorithm \eqref{pushsumx}--\eqref{pushsumy}. 

The following results show that important properties of the standard push-sum still hold for the generalized algorithm \eqref{eq:convexx}--\eqref{eq:convexy}, with variant appearance. 

\begin{proposition}\label{prop:convex}
    If $\{ \bbb{G}(t) \}$ is uniformly strongly connected, then for any fixed $\tau\ge0$, $S(t)\cdots S(\tau+1)S(\tau)$ will converge to $ \1 y^\top(\tau)$ exponentially fast as $t\rightarrow\infty$. 
\end{proposition}

The proposition can be proved using similar arguments to those in the proof of Proposition~\ref{lemma:expconvergen}.

{\bf Proof of Proposition~\ref{prop:convex}:}
Since \eqref{eq:convexy} shares the same $y_i$ update as \eqref{pushsumy} except for the initial value setting, the result of Lemma \ref{lemma:yixuan} still holds for \eqref{eq:convexy}, and the result of Lemma \ref{lemma:y_bound} can be further polished as $1 \ge y_i(t) \ge c_i\eta$ for all $i$ and $t$. 
From \eqref{mu}, there exist constants $c>0$ and $\mu \in [0,1)$ and a sequence of stochastic vectors $ \{ v(t)\}$ such that 
$
    | [\Phi_W(t+1,\tau)]_{ij} - v(t) |\le c \mu^{t-\tau}
$ for all $i,j \in \scr V$ and $t \ge \tau \ge 0$.
Note that $\sum_{i=1}^n y_i(t)$ always equals $\sum_{i=1}^n c_i = 1$ and, by Lemma~\ref{lemma:y_bound},  all $y_i(t)$ are always positive. From Lemma~\ref{lemma:yixuan}, for all $t\ge\tau\ge0$,
\begin{align*}
&\Big| [\Phi_S(t+1,\tau)]_{ij} - y_j(\tau) \Big|
= \Big| \frac{ y_j(\tau) [\Phi_W(t+1,\tau) ]_{ij} }{y_i(t+1)}  - y_j(\tau) \Big|\\
=\;&\Big| \frac{ y_j(\tau) [\Phi_W(t+1,\tau) ]_{ij} - y_j(\tau) [\Phi_W(t+1,\tau)y(\tau)]_i}{y_i(t+1)}   \Big|
\\
=\;& \Big| \frac{ y_j(\tau) \left([\Phi_W(t+1,\tau)]_{ij} - v_i(t) + v_i(t) \right) - y_j(\tau)\sum_{k=1}^n ([\Phi_W(t+1,\tau)]_{ik} - v_i(t) + v_i(t) ) y_k(\tau) }{y_i(t+1)} \Big|\\
=\;&  \Big| \frac{  y_j(\tau) ([\Phi_W(t+1,\tau)]_{ij} - v_i(t) ) - y_j(\tau) \sum_{k=1}^n ( [\Phi_W(t+1,\tau)]_{ik}- v_i(t) ) y_k(\tau)    }{ y_i(t+1) } \Big|\\
\le \;&   \frac{   y_j(\tau) \left|[\Phi_W(t+1,\tau)]_{ij} - v_i(t) \right| +  y_j(\tau) \sum_{k=1}^n \left| [\Phi_W(t+1,\tau)]_{ik}- v_i(t) \right| y_k(\tau) }{ y_i(t+1) } \\
\le \;&  \frac{  y_j(\tau) c \mu^{t-\tau} + y_j(\tau) \sum_{k=1}^n c \mu^{t-\tau} y_k(\tau)     }{ y_i(t+1) } 
= \frac{ 2 y_j(\tau) c \mu^{t-\tau} }{ y_i(t+1) } \le  \frac{2c}{ c_i\eta } \mu^{t-\tau} \le  \frac{2c}{ \eta \min_{i\in\scr V} c_i} \mu^{t-\tau},
\end{align*}
which immediately implies the proposition.
\hfill$\qed$

Proposition~\ref{prop:convex} immediately implies, by setting $\tau=0$, that 
$S(t)\cdots S(1)S(0)$
    will converge to $\1 [c_1 \; c_2\; \cdots \; c_n]$ exponentially fast as $t\rightarrow\infty$.
This leads to Theorem \ref{thm:convex}.

{\bf Proof of Theorem \ref{thm:convex}:}
Note that the equation \eqref{eq:update_ratio} still holds for the algorithm \eqref{eq:convexx}--\eqref{eq:convexy}. Then, 
$z(t+1)=S(t)z(t)=S(t)\cdots S(1)S(0)z(0)$. 
From the preceding, it is easy to show that $z(t+1)$ will converge to $\1 [c_1 \; c_2\; \cdots \; c_n] z(0) = (\sum_{k=1}^n c_k x_k^{{\rm int}})\1$ exponentially fast as $t\rightarrow\infty$. 
\hfill$\qed$

More can be said. The key property of the algorithm \eqref{eq:convexx}--\eqref{eq:convexy} is as follows.

\begin{proposition} \label{piconvex}
    If $\{ \bbb{G}(t) \}$ is uniformly strongly connected, then the sequence of stochastic matrices $\{S(t)\}$ has a unique absolute probability sequence $\{\pi(t)\}$ with
    $\pi_i(t)=y_i(t)$ for all $i\in\scr V$ and $t \ge 0$. 
\end{proposition}

Proposition~\ref{piconvex} can be proved using similar arguments to those in the proof of Proposition~\ref{lemma:push-sum_pi_intfty}.

{\bf Proof of Proposition~\ref{piconvex}:}
From Proposition~\ref{prop:convex},  $\{S(t)\}$ is ergodic, and thus has a unique absolute probability sequence $\{\pi(t)\}$. 
Since $\{y(t)\}$ is a sequence of stochastic vectors, to prove the proposition is equivalent to show that $y^\top(t) = y^\top(t+1) S(t)$ for all $t\ge 0$. 
To see this, from \eqref{eq:s} and Assumption~\ref{assum:weighted matrix}, 
$
    [y^\top(t+1) S(t)]_j 
    = \sum_{i=1}^n y_i(t+1) s_{ij}(t)  
    = \sum_{i=1}^n w_{ij}(t)y_j(t)
    = y_j(t) 
$
for all $j \in \scr V$. 
\hfill$\qed$


Proposition~\ref{piconvex} can be viewed as a generalization of Proposition~\ref{lemma:push-sum_pi_intfty}, the key property of the standard push-sum. To see this, consider the algorithm \eqref{eq:convexx}--\eqref{eq:convexy} without assuming $\sum_{i=1}^n c_i = 1$, that is, all $c_i$ are all positive but their sum $\kappa \dfb\sum_{i=1}^n c_i$ may not equal one. 
Using similar arguments to those in the proofs of Propositions~\ref{prop:convex} and \ref{piconvex}, we can show that if $\{ \bbb{G}(t) \}$ is uniformly strongly connected, then for any fixed $\tau\ge0$, $S(t)\cdots S(\tau+1)S(\tau)$ will converge to $ \frac{1}{\kappa}\1 y^\top(\tau)$ exponentially fast as $t\rightarrow\infty$, and more importantly, the sequence of stochastic matrices $\{S(t)\}$ has a unique absolute probability sequence $\{\pi(t)\}$ with
    $\pi_i(t)=\frac{y_i(t)}{\kappa}$ for all $i\in\scr V$ and $t \ge 0$. 
Consequently, each $x_i(t)/y_i(t)$ converges to $\frac{1}{\kappa}\sum_{k=1}^n c_k x_k^{{\rm int}}$ exponentially fast as $t\rightarrow\infty$. 
From this point of view, the results of Proposition~\ref{lemma:push-sum_pi_intfty} and Proposition~\ref{piconvex} are consistent, while the latter corresponds to a more general push-sum algorithm.


The above properties guarantee that the modified $y$ update
\eqref{eq:convexy} can be applied to push-sum based distributed optimization algorithms with the same analysis tool presented in the current paper. 
For example, by setting $y_i(0)=c_i$ instead of 1, the subgradient-push algorithm \eqref{eq:pushsub_x}--\eqref{eq:pushsub_y} becomes
\begin{align*}
    x_i(t+1) &= \sum_{j\in\scr{N}_i(t)} w_{ij}(t)\Big[x_j(t) - \alpha(t)g_j(t)\Big], \;\;\;\;\; x_i(0)\in\R^d,\label{eq:pushsub_x}\\  y_i(t+1) &= \sum_{j\in\scr{N}_i(t)} w_{ij}(t)y_j(t),\;\;\;\;\; y_i(0)=c_i,
\end{align*}
which, with $\sum_{i=1}^n c_i = 1$, can also solve the distributed convex optimization problem $f(z)=\frac{1}{n}\sum_{i=1}^n f_i(z)$, using the same arguments as in the analysis of the subgradient-push algorithm \eqref{eq:pushsub_x}--\eqref{eq:pushsub_y}. The same modification also works for stochastic gradient-push and heterogeneous subgradient.

The preceding discussion implies that changing initial values of subgradient-push only cannot solve the distributed weighted optimization $\sum_{i=1}^n c_if_i(z)$ with all $c_i$ being positive. In fact, the distributed weighted optimization problem can be straightforwardly solved by assuming each agent $i$ knows $c_i$ and then treating $c_if_i(z)$ as the private convex cost function of agent $i$. 




\section{Concluding Remarks}

Push-sum is a well-known algorithm for reaching an average consensus over time-varying directed graphs. For general discrete-time linear consensus processes, it is known that a consensus will be reached exponentially fast if and only if the underlying neighbor graph sequence is ``repeatedly jointly rooted'' \cite{cdc14}, a weaker connectivity condition compared with ``repeatedly jointly strongly connected'' or uniformly strongly connected. Olshevsky and Tsitsiklis proved that there exists no common quadratic Lyapunov function for such time-varying linear consensus systems \cite{noquadratic}. 
A few years later, Touri and Nedi\'c constructed a time-varying quadratic Lyapunov comparison function, based on the concept of an absolute probability sequence of a sequence of stochastic matrices (cf. Definition \ref{def: absolute prob}), for discrete-time linear consensus processes over uniformly strongly connected graphs \cite{touri2012product,touritac,tacrate}.
Although (uniformly) strong connectedness is not broad enough for general consensus, it is a necessary and sufficient condition for distributed averaging and distributed optimization. To our knowledge, there is no existing quadratic Lyapunov (comparison) function for analyzing the push-sum algorithm and push-sum based distributed algorithms. 
Almost all existing push-sum based distributed optimization algorithms more or less rely on the analyses in the pioneering work \cite{nedic} by Nedi\'c and Olshevsky. 

In this paper, we have revisited the push-sum algorithm and provided an alternative convergence proof by characterizing the dynamics of all $x_i/y_i$ ratios as a nonlinear consensus process and establishing its unique absolute probability sequence with an explicit expression. This explicit absolute probability sequence has then been used to construct quadratic Lyapunov (comparison) functions for push-sum based distributed first-order optimization algorithms. Appealing to this new analysis tool, we have also revisited subgradient-push and stochastic gradient-push, two important algorithms for distributed convex optimization over unbalanced directed
graphs, and proved that they both converge at an optimal rate even over time-varying graphs, which improves the existing results by an order of $O(\ln t)$. These two algorithms illustrate that push-sum based distributed optimization can achieve the best possible convergence rate for time-varying, unbalanced, directed graphs, without any network-wide information. 
As one future direction, the proposed tool is expected to be applicable to analyze other push-sum based optimization algorithms and improve/simplify their convergence analyses, for example, DEXTRA \cite{xi2017dextra} and Push-DIGing \cite{nedic2017achieving}. 

The proposed novel analysis tool has further assisted us to design a heterogeneous distributed subgradient algorithm with a unified convergence analysis, subsuming all existing push-sum and average consensus based subgradient algorithms as special cases. The heterogeneous algorithm allows each agent to arbitrarily switch between subgradient-push and push-subgradient at any time, which is expected to be beneficial to protect privacy against an honest-but-curious adversary or an external eavesdropping adversary.
The tool has also been extended to analyze distributed weighted averaging, 
subsuming distributed averaging as a special case.

While being applied only to distributed optimization, the proposed tool is expected to be useful to analyze other push-sum based distributed algorithms. In our recent work \cite{finite}, the developed novel push-sum property (cf. Proposition \ref{lemma:push-sum_pi_intfty}) has been used to establish the boundedness of a push-sum based distributed linear stochastic
approximation process, which can be directly applied to distributed temporal difference (TD) learning over time-varying unbalanced directed graphs.
Future directions include extending the proposed tool to analyze other push-sum based algorithms (e.g., non-convex optimization, reinforcement learning, deep learning) and to deal with more realistic scenarios (e.g., communication delays, asynchronous updating, package drops). In particular, it is worth noting that a converge rate gap of order $O(\ln t)$ also exists for asynchronous push-sum based distributed gradient \cite{asyngradpush}, and we expect the novel tool here can be tailored to fill this gap.



\section*{Acknowledgement}

The authors wish to thank Jie Lu (ShanghaiTech University) Wenhan Gao (Stony Brook University), and Jinglin Yang (Stony Brook University) for useful discussion.

\bibliographystyle{unsrt}
\bibliography{push}
\end{document}